\documentclass[a4paper,10pt]{article}
\usepackage[english, french]{babel}
\usepackage{amsfonts, amsmath,amscd,amssymb} 
\usepackage{latexsym}
\usepackage[latin1]{inputenc}
\usepackage{url}
\makeatletter
\@addtoreset{equation}{section}

\makeatother
 \usepackage[bookmarks]{hyperref}

\newcommand{\C}{\mathbb{C}}
\newcommand{\R}{\mathbb{R}}
\newcommand{\N}{\mathbb{N}}
\newcommand{\Z}{\mathbb{Z}}

\newcommand{\im}{\mathrm{Im}\,}

\begin{document}
\title{Loi de Weyl presque sûre pour un système différentiel en dimension 1 }
\author{William Bordeaux Montrieux \\ 
\small Centre de Mathématiques Laurent Schwartz\\ 
\small Ecole Polytechnique\\
\small FR 91120 Palaiseau cedex\\ 
\small bordeaux@math.polytechnique.fr}
\date{}
\maketitle

\newtheorem{theo}{Th\'{e}orème}[section]
\newtheorem{corol}[theo]{Corollaire}
\newtheorem{prop}[theo]{Proposition}
\newtheorem{lemme}[theo]{Lemme}
\newtheorem{definition}[theo]{D\'{e}finition}
\newtheorem{hypothese}[theo]{Hypothèse}
\newtheorem{remarque}[theo]{Remarque}

\begin{abstract} 
Nous considérons une classe assez générale de  systèmes différentiels sur le cercle
avec une perturbation aléatoire d'ordre inférieur. Nous adoptons deux points de vue,
semiclassique et haute fréquence.
Nous montrons (a) que dans la limite $h\to 0,$ les valeurs propres se distribuent selon 
une loi de Weyl avec une probabilité très proche de 1, (b) que les grandes valeurs propres
se distribuent \textit{presque sûrement}  selon une loi de Weyl. 
\end{abstract}
\selectlanguage{english}   
\begin{abstract} 
We consider quite general differential operators on the circle with a small random lower order
perturbation. We embrace two points a view, the semiclassical and the high energy limits. We show
(a) in the semiclassical limit, that the eigenvalues inside a subdomain of the pseudospectrum are distributed according to a Weyl law with a probability close to 1, (b) that the large eigenvalues  obey 
 a Weyl law \textit{almost surely}.
\end{abstract}
\selectlanguage{french}  
\tableofcontents
\section{Introduction}
Les constructions de quasimode  de E.B.~Davies \cite{Da}, M.~Zworski \cite{Zw1}
et d'autres \cite{DeSjZw, Pr} impliquent 
que les opérateurs $h$-pseudodifférentiels non-auto\-adjoints ont, en général,
la norme de la résolvante qui est très grande
lorsque le paramètre spectral
$z$ se déplace à l'intérieur de l'image du symbole principal. Dit autrement,
le spectre est très instable sous petites perturbations.
Une question naturelle est de comprendre comment les valeurs propres
bougent quand l'opérateur est perturbé, et notamment lorsque la
perturbation est aléatoire.

Dans \cite{Ha2}, M. Hager  considère certaines classes d'opérateurs 
pseudodifféren\-tiels semiclassiques $P$ sur $\R$; incluant les opérateurs 
différentiels. Elle utilise des petites perturbations aléatoires multiplicatives 
$\delta Q_{\omega},$ où $\delta$ est un petit paramètre. Soit un domaine $\Gamma\Subset\C$
avec une frontière lisse, on suppose que $p^{-1}(z)$ est une collection finie 
de points pour $z$ dans $\Gamma$ et pour lequel
$\{p,\bar{p}\}(\rho)\ne 0$ si $\rho\in p^{-1}(\Gamma).$
Sous des hypothèses additionnelles, Hager a montré 
qu'avec une probabilité qui tend vers 1 lorsque $h\to 0,$
pour $\delta=e^{-\epsilon/h},$ les valeurs propres
de l'opérateur perturbé se distribuent selon une Loi de Weyl dans $\Gamma,$
\textit{ce qui était déjà bien connu dans le cas autoadjoint,} 
\begin{equation}\label{ii.2}
|\#(\sigma(P+\delta Q_{\omega})\cap\Gamma)-
\frac{1}{2\pi h}\mathrm{vol}\, (p^{-1}(\Gamma))|\le \frac{C\sqrt{\epsilon}}{h},
\quad h\to 0.
\end{equation}

Mentionnons que  M.~Hager, J.~Sj\"ostrand  ont étendu ce résultat au cas
des opérateurs sur $\R^{n}$ \cite{HaSj}, et que J.~Sj\"ostrand l'a lui étendu au cas des variétés
compactes \cite{Sj}.

Dans ce travail, nous allons étudier des systèmes elliptiques d'opérateurs 
différentiels sur $S^{1}$ avec des perturbations aléatoires. En adaptant des techniques
d'Hager \cite{Ha2} nous allons d'abord établir une loi de Weyl avec une probabilité 
proche de 1 dans le cas semiclassique avec des petites perturbations.
Ensuite dans le cas non-semiclassique ($h=1$) nous montrerons que les 
grandes valeurs propres se distribuent \textit{presque sûrement} selon la loi de Weyl.\\

\noindent\textbf{Remerciements.} 
Ce travail fait partie de la thèse préparée sous la direction de J.~Sj\"ostrand.
L'auteur tient aussi à remercier deux personnes M.~Zworski pour son accueil et les discussions
très utiles lors de son séjour à Berkeley, et le rapporteur des Annales de l'institut Henri Poincaré
qui est pour beaucoup dans cette nouvelle version.

\section{Enoncé des résultats}
\paragraph{Asymptotique semiclassique.}
Considérons l'opérateur différentiel non-auto\-adjoint dans $L^{2}(S^{1},\C^{n})$
\begin{equation}\label{i.1}
P(h)=\sum_{0\le\alpha\le m}A_{\alpha}(x;h)(hD_{x})^{\alpha},\quad 
h\in (0,1],\quad D_{x}=\frac{1}{i}\frac{\partial}{\partial x},
\end{equation}
où chaque $A_{\alpha}$ est une matrice $n\times n$ complexe dépendant de manière 
$C^{\infty}$ de $x$, et
admettant la repésentation asymptotique dans $C^{\infty}(S^{1})$,
\begin{equation}\label{i.2}
A_{\alpha}(x;h)\sim A_{\alpha,0}(x)+hA_{\alpha,1}(x)+h^{2}A_{\alpha,2}(x)+\ldots,\quad h\to 0.
\end{equation}
Le domaine de définition $\mathcal{D}(P)$ choisi pour $P$ est l'espace de Sobolev semiclassique
$H_{sc}^{m}(S^{1},\C^{n})$ défini par 
\begin{equation}\label{i.3}
\left\{u\in L^{2}(S^{1},\C^{n})\bigg|\,\|u\|_{m,h}^{2}=\sum_{0\le\alpha\le m}\|(hD_{x})^{\alpha}u\|^{2}
<\infty\right\}.
\end{equation}
Le symbole principal semiclassique de $P$ est donné par 
\begin{equation}\label{i.4}
p(x,\xi):=
\sum_{0\le\alpha\le m} A_{\alpha,0}(x)\xi^{\alpha},\quad (x,\xi)\in T^{\ast}S^{1}.
\end{equation}
\begin{hypothese}\label{i1}
On suppose que $P$ est elliptique (au sens où $\det A_{m,0}$ ne s'annule pas).
\end{hypothese}

Nous notons l'ensemble des valeurs propres du symbole principal $p$ par
\begin{equation}
\Sigma(p) =\bigcup_{(x,\xi)\in T^{\ast}S^{1}}\sigma(p(x,\xi)),
\quad \sigma(p(x,\xi)):=\mbox{spectre de }p(x,\xi).
\label{i.5} 
\end{equation}
Si $P$ est un opérateur scalaire, alors $\Sigma$ est l'ensemble des valeurs de $p$.
\begin{prop}\label{i2}
Sous l'hypothèse précédente, pour tous $z$, $P-z:\mathcal{D}(P)\to L^{2}(S^{1})$ 
est un opérateur de Fredholm d'indice zéro. 
\end{prop}
\textbf{Preuve.}
Il est connu qu'un opérateur elliptique sur une variété compacte 
(ici $S^{1}$) est de Fredholm.  Après multiplication par 
$A_{m}^{-1},$ on est ramené au cas où $A_{m}=I.$ 
Puis, en utilisant l'invariance de l'indice de Fredholm par déformation elliptique,
on obtient que l'indice de $P-z$ est égal à celui de $(hD)^{m};$
les termes de degré inférieur ont été écrasés.
Pour finir, il est clair que l'indice de $(hD)^{m}$ est zéro.
\hfill $\square$ \medskip

En particulier, s'il existe un point $z_{0}$ pour lequel la résolvante $(P-z_{0})^{-1}$
existe (ce qui est toujours le cas si $\Sigma(p)\ne \C$), 
alors nous trouvons que le spectre est discret dans $\C.$  
En effet, par la théorie de Fredholm analytique, nous savons que 
pour un  opérateur $A$ d'indice zéro, dont le spectre n'est pas égal à $\C,$  
alors le spectre consiste en des valeurs propres discrètes. \\

Pour $z$ fixé, $q_{z}(x,\xi)$  désigne, dans la suite, le déterminant de $p(x,\xi)-z.$ 
Nous définissons l'ensemble
\begin{equation}
\Phi =\{z\in \Sigma\,|\,\exists(x,\xi)\in T^{\ast}S^{1}\mbox{ avec } z\in\sigma(p(x,\xi))
\mbox{ et } \{q_{z},\bar{q}_{z}\}(x,\xi)=0 \} \label{i.6}
\end{equation}
où $\{\bullet,\bullet\}$ désigne le crochet de Poisson. 
$\Sigma,\Phi$ sont fermés 
et   $\Lambda(p):=\Sigma\setminus\Phi$ 
est un ensemble ouvert. 

Nous montrerons, dans la proposition \ref{sp4},  que l'image réciproque de 
zéro par $q_{z}$ pour 
$z$ donné dans $\Lambda(p),$ est un 
ensemble de la forme
\begin{equation}\label{i.7}
\forall z\in \Lambda(p),\quad
q^{-1}_{z}(0)
=\{\rho_{+}^{\nu}(z), \rho_{-}^{\nu}(z)|\,\nu=1,\dots,\beta(z)\},
\end{equation}
où $\beta(z)<\infty$ est  localement constant, et
\begin{equation}\label{i.8}
\pm \frac{1}{2i}\{q_{z},\bar{q}_{z}\}(\rho_{\pm})>0.
\end{equation}

Ce qui implique que pour tout $z\in\Lambda(p),$ $\nu=1,\ldots,\beta,$
\[\exists\, e_{+}^{\nu}=e_{+}^{\nu}(x,z;h)\in \mathcal{S},\,
\|e_{+}^{\nu}\|=1, \,\|(P-z)e_{+}^{\nu}\|=\mathcal{O}(h^{\infty}),\]
$e_{+}^{\nu}$ est une solution BKW concentrée près de $\rho_{+}^{\nu},$
et 
\[\exists\, e_{-}^{\nu}=e_{+}^{\nu}(x,z;h)\in \mathcal{S},\,
\|e_{-}^{\nu}\|=1, \,\|(P-z)^{\ast}e_{+}^{\nu}\|=\mathcal{O}(h^{\infty}),\]
$e_{-}^{\nu}$ est une solution BKW concentrée près de $\rho_{-}^{\nu}.$

\begin{hypothese}\label{i3} 
Soit $\Omega\Subset\Lambda(p)$ et connexe. 
On demande que pour tout $z\in\Omega,$
\begin{equation}
\rho_{\pm}^{\nu}(z)=(x^{\nu}(z),\xi_{\pm}^{\nu}(z)),\quad
x^{\nu}\ne x^{\kappa},\,\nu\ne\kappa.\label{i.9}
\end{equation}
et que $\xi_{+}^{\nu}\ne 0$ pour tout $\nu\in1,\ldots,\beta$,
où $\beta$ est la valeur constante de $\beta(z)$ sur la composante
connexe de $\Lambda(p)$ contenant $\Omega$.
\end{hypothese}

Soient  $(\mathcal{M},\mathcal{A},\mathbb{P})$ un espace de probabilité et
$Q_{\omega}$ est  un  opérateur différentiel d'ordre inférieur à $m$  de 
$L^{2} (S^{1} )$ dans lui-même, de domaine dense,
\begin{equation}\label{i.11}
Q_{\omega}=\sum_{\alpha_{0}\le\alpha\le\alpha_{1}}Q_{\alpha}(x;h)(hD_{x})^{\alpha},\quad
0\le\alpha_{0}\le\alpha_{1}\le m-1.
\end {equation} 
Ici $(Q_{\alpha}^{i,j})_{i,j}$ est une matrice $n\times n$ où chaque élément 
est une série de Fourier aléatoire,
c'est-à-dire 
\begin{equation}\label{i.12}
Q_{\alpha}^{i,j}(x;h)=\sum_{k\in\Z}q_{\alpha,k}^{i,j}(h)\frac{e^{ikx}}{\sqrt{2\pi}}.
\end{equation}
\begin{remarque}
De manière général, nous adoptons la convention  
suivante : les coefficients d'une matrice 
$Q$ seront indiqués par les exposants $i, j$, $Q^{i,j}.$ 
\end{remarque}

On adopte l'hypothèse suivante sur les variables aléatoires $q_{\alpha,k}^{i,j}$:
\begin{hypothese}\label{i4}
Les coefficients de Fourier $q^{i,j}_{\alpha,k}$ sont des variables aléatoires 
(pour faire court v.a.) complexes indépendantes de  loi 
$\mathcal{N}(0,(\sigma_{\alpha,k}^{i,j})^{2}).$
La variance peut dépendre de $h.$ Pour tout $i,j, \alpha$ et $0<h\le 1,$
\begin{equation}\label{i.13}
\sigma^{i,j}_{\alpha,k}(h)\le
\tilde{C}\langle k \rangle^{-\rho},
\end{equation}
et pour $\alpha=\alpha_{1},$ nous avons pour tout $i,j$
\begin{equation}\label{i.14}
 \sigma^{i,j}_{\alpha_{1},k}(h)\ge
 \frac{1}{\tilde{C}}\langle k\rangle^{-\rho},
 \end{equation}
où les constantes $\tilde{C}>0$ et $\rho>1$ sont indépendantes de $\alpha,i,j,k,$ et $h$
et où on utilise la notation standard $\langle k\rangle=(1+|k|^{2})^{\frac{1}{2}}.$
\end{hypothese}

On rappelle  que $X$ suit une loi gaussienne $\mathcal{N}(m,\sigma^{2})$
complexe d'espérance $m\in \C$ et de variance $\sigma^{2}>0,$
si sa densité est 
\[\varphi(z)=
\left\{\begin{array}{l}
\frac{1}{\pi \sigma^{2}}e^{\frac{-|z-m|^{2}}{\sigma^{2}}},\quad \sigma>0,\\
\delta(z-m)\,(\mbox{masse de Dirac en }z=m), \quad \sigma=0.
\end{array}
\right.\]
La propriété remarquable des v.a. gaussiennes est que 
la somme de deux v.a. gaussiennes indépendantes reste une v.a. gaussienne où
les espérances et les variances s'additionnent respectivement. 

Sous ces conditions, $Q_{\omega}$ est presque sûrement (p.s.)  borné comme 
opérateur de $H^{m}_{sc}$ dans $L^{2}.$
Ce fait découle du résultat suivant concernant la régularité des fonctions $Q_{\alpha}^{i,j}(x;h)$ :
\begin{prop}\label{i5} Sous l'hypothèse précédente, pour chaque $\alpha,i,j,$ 
$Q^{i,j}_{\alpha}(x;h)$ 
repré\-sente p.s. une fonction continue.
\end{prop}
\textbf{Preuve.} Il suffit de remarquer, grâce à l'inégalité de Markov, que 
\[\mathbb{P}(\sum_{k\in\Z}|q_{\alpha,k}^{i,j}|> t)\le t^{-1}
\mathbb{E}(|X|)\sum_{k\in \Z}\sigma_{\alpha,k}^{i,j},\]
où $X$ suit une loi gaussienne standard $\mathcal{N}(0,1).$  On fait ensuite tendre $t$ 
vers l'infini pour voir que la série  aléatoire (\ref{i.9}) converge normalement presque sûrement,
d'où la continuité.
\hfill $\square$ \medskip 

Ils existent des résultats très fins concernant la régularité, l'irrégularité des 
séries de Fourier aléatoires gaussiennes, voir \cite{Ka}. 

Nous introduisons, pour $(x,\xi)\in T^{\ast}S^{1}$ et $\Gamma\Subset\C$ donné,
le nombre de valeurs propres de $p(x,\xi)$ dans $\Gamma$ par
\begin{equation}
m_{\Gamma}(x,\xi):=\#(\sigma(p(x,\xi))\cap \Gamma).
\end{equation}
Nous nous proposons alors d'établir le résultat suivant :
\begin{theo}\label{i6}
Supposons admis les hypothèses  \ref{i1}, \ref{i3} et \ref{i4}. 
Soit $\Gamma\Subset\Omega$ un ouvert 
à bord $C^{2}$ par morceaux. 
On entend par cela que $\partial\Gamma$ peut être paramétré 
par une courbe $S^{1}$ dans $\C$ continue et $C^{2}$ en dehors d'un 
nombre fini de points $a_{1},a_{2},\ldots,$
et pour lesquels l'angle formé par la dérivées à gauche et à droite est non nul. 
Pour tous $\gamma_{1},N_{0}>0,$ il existe une 
constante positive $C>0$ telle que, pour
$h^{N_{0}}<\delta<h^{\rho+\gamma_{1}+\frac{1}{2}}|\ln h|^{-2},$
le spectre de $P-\delta Q_{\omega}$ est discret et on a
\[
|N(P-\delta Q_{\omega},\Gamma)-\frac{1}{2\pi h} \iint m_{\Gamma}(x,\xi)\,dxd\xi |\le C
h^{-\frac{1}{2}}|\ln h|^{\frac{1}{2}},
\]
avec une probabilité 
\[\ge1-Ch^{2\gamma_{1}}|\ln h|^{-\frac{1}{2}}.\]
\end{theo}

Notons que lorsque $\alpha_{1}=\alpha_{0}=0,$ nous nous trouvons 
dans la situation d'une perturbation multiplicative aléatoire.
Nous donnerons  au  théorème  \ref{ss9} une version de la loi de Weyl
pour une famille $\mathcal{G}$  de domaine $\Gamma$ dans $\Omega.$

\paragraph{Asymptotique des grandes valeurs propres.}
Soit l'opérateur différentiel non-autoadjoint dans $L^{2}(S^{1},\C^{n})$
\begin{equation}\label{i.17}
P=\sum_{0\le\alpha\le m}^{}A_{\alpha}(x)D_{x}^{\alpha},
\quad A_{\alpha}(x)\in C^{\infty}(S^{1}).
\end{equation}

Le domaine de définition naturel est l'espace de Sobolev $H^{m}(S^{1},\C^{n}).$
On impose comme précédemment une hypothèse d'ellipticité
\begin{equation}\label{i.18}
\det A_{m}(x)\ne 0,\quad x\in S^{1},
\end{equation}
rendant l'opérateur $P-z$ de Fredholm d'indice zéro pour tout $z.$
En particulier, si $P-z$ est bijectif pour au moins une valeur de $z,$
et nous trouvons donc  que le  spectre de $P$ est discret. 
Le symbole principal classique de $P$ est
$p_{m}(x,\xi):=A_{m}(x)\xi^{m},$
et nous désignons par $\Sigma(p_{m})$ l'ensemble 
des valeurs propres de $p_{m},$ c'est à dire
\begin{equation}\label{i.20}
\Sigma(p_{m})=\bigcup_{(x,\xi)\in T^{\ast}S^{1}}\sigma(p_{m}(x,\xi)).
\end{equation}

Pour $z$ donné, on écrit  $q_{m,z}(x,\xi)$ pour $\det (p_{m}(x,\xi)-z).$ 
Nous introduisons ensuite 
l'ensemble,
\begin{equation}\label{i.21}
\Phi=\{z\in\Sigma|\,\exists(x,\xi)\in T^{\ast}S^{1}\mbox{ avec }
z\in\sigma(p_{m}(x,\xi))\mbox{ et }\{q_{m,z},\overline{q_{m,z}}\}(x,\xi)=0\}.
\end{equation}

Nous utilisons la perturbation,
\begin{equation}\label{i.22}
Q_{\omega}=\sum_{\alpha_{0}\le\alpha\le \alpha_{1}}Q_{\alpha}(x)
D_{x}^{\alpha},\quad
0\le\alpha_{0}\le\alpha_{1}\le m-1,
\end{equation}
où chaque élément $Q^{i,j}_{\alpha}$ est une série de Fourier aléatoire
\[Q_{\alpha}(x)=\sum_{k\in\Z}q_{\alpha,k}^{i,j}(x)\frac{e^{ikx}}{\sqrt{2\pi}}.\]
Nous supposons de plus que les coefficients $q^{i,j}_{\alpha,k}$ vérifient 
l'hypothèse \ref{i4}.
La proposition \ref{i5} nous dit alors que $Q_{\omega}$ est un opérateur différentiel 
dont les coefficients sont p.s. continus. De plus,
 puisque $P$ et $P-Q_{w}$ ont le même symbole principal, alors  p.s.
$P-Q_{w}$ est un opérateur de Fredholm d'indice zéro.\\

Nous sommes  intéressés ici par la distribution des grandes valeurs propres de $P-Q_{\omega}$
dans les dilatés d'un profil conique inclus dans $\Lambda(p_{m}):=\Sigma\setminus\Phi$
(qui est un cône du fait de l'homégénéité du symbole principal).
Choisissons, $\Omega,$ un cône ouvert connexe dans $\Lambda(p_{m}).$ 

Pour $z$ fixé dans $\Lambda(p_{m}),$
l'image réciproque de zéro par $q_{m,z}$ est un ensemble de la forme
\begin{equation}\label{i.23}
\forall z\in \Lambda(p_{m}),\quad
q^{-1}_{m,z}(0)=\{\rho_{+}^{\nu}(z), \rho_{-}^{\nu}(z)|\,\nu=1,\dots,\beta(z)\},
\end{equation}
où $\beta(z)<\infty$ est  constant sur chaque composante connexe de $\Lambda(p_{m})$, et
\begin{equation}\label{i.24}
\pm \frac{1}{2i}\{q_{m,z},\bar{q}_{m,z}\}(\rho_{\pm})>0.
\end{equation}  
On fait alors l'hypothèse suivante :
\begin{hypothese}\label{i10} 
On demande que pour tout $z\in\Omega,$
\begin{equation}
\rho_{\pm}^{\nu}(z)=(x^{\nu}(z),\xi_{\pm}^{\nu}(z)),\quad
x^{\nu}\ne x^{\kappa},\,\nu\ne\kappa.\label{i.25}
\end{equation}
\end{hypothese}
Puisque le symbole principal est homogène par rapport à $\xi$,
nous avons forcément $\xi_{+}^{\nu}(z)\ne 0$ si $z\ne 0$.

Soient $\theta_{1}^{0}$ et $\theta_{2}^{0}$ tels que  
\[ \Lambda(p_{m})\supset \Omega=\{ re^{i\theta}|\,
r>0,\, \theta_{1}^{0}<\theta<\theta_{2}^{0}\}.\]
Prenons $\theta_{1},\theta_{2}\in ]\theta_{1}^{0},\theta_{2}^{0}[,$
avec $\theta_{1}\le\theta_{2},$ et
$g,h\in C^{2}([\theta_{1},\theta_{2}],\R_{+})$ satisfaisant $h<g.$ 
Nous introduisons alors l'ensemble
\begin{equation}\label{i.26}
\Omega\Supset\Gamma_{\theta_{1},\theta_{2}}(h,g):=
\{re^{i\theta}|\,\theta_{1}\le \theta\le\theta_{2},\,h(\theta)\le r\le g(\theta)\}.
\end{equation}
Pour $\theta_{1},\theta_{2}$ fixés on écrira parfois
$\Gamma(h,g)$ à la place $\Gamma_{\theta_{1},\theta_{2}}(h,g).$
Nous notons pour tout $(x,\xi)\in T^{\ast}S^{1}$ et $\Gamma\subset\C$
\begin{equation}
m_{\Gamma}(x,\xi):=\#(\sigma(p_{m}(x,\xi))\cap \Gamma).
\end{equation}

Notre résultat principal est le suivant :
\begin{theo}\label{j8}
Soit $\Omega$ un cône connexe dans $\Lambda(p_{m}).$ On suppose que  l'hypothèse
d'ellipticité est satisfaite et  \ref{i4}, \ref{i10} sont vérifiées. Si $m-\alpha_{1}-\rho-\frac{3}{4}>0,$
alors il existe $\widetilde{C}>0$ et
$\widetilde{M}\subset\mathcal{M}$ avec $\mathbb{P}(\widetilde{M})= 1$ tels que,
pour tout   $\omega\in \widetilde{M}$,
le spectre de $P-Q_{\omega}$ est discret, et 
le nombre $N(P-Q_{\omega},\lambda\Gamma(0,g))$
de valeurs propres de $P-Q_{\omega}$ dans
 $\lambda\Gamma(0,g)\Subset\Omega$ satisfait
\begin{align*}
&\forall \lambda\ge 0,\\
&|N(P-Q_{\omega},\lambda\Gamma(0,g))
-\frac{1}{2\pi}
\iint m_{\lambda\Gamma(0,g)}(x,\xi)\,dxd\xi|
\le C(\omega)+\widetilde{C} \lambda^{1/(2m)}\sqrt{\ln \lambda}. 
\end{align*}
La constante $C(\omega)<+\infty$ dépend de $\omega$, mais pas de 
$\lambda$. 
\end{theo}

Pour le cas $m=2,$ le théorème n'est vérifié que pour des 
perturbations multiplicatives avec $1<\rho<\frac{5}{4}.$\\

Notre démonstration est organisée comme suit. 
Après quelques rappels sur les opérateurs pseudodifférentiels
et les notations utilisées, nous montrons que, 
pour tout $z\in\Sigma\setminus\Phi,$
l'ensemble $q^{-1}_{z}(0)$ est composé par autant de points 
$\rho_{+}^{\nu}$ et $\rho_{-}^{\nu}$ tels que le crochet de Poisson 
$\{q_{z},\overline{q_{z}}\}$ soit strictement positif aux points $\rho_{+}^{\nu}$
et strictement négatif aux points $\rho_{-}^{\nu}$ (section 3).
Dans \cite{Ha2}, Hager fait l'hypothèse qu'il y a autant de point $\rho_{+}$
que $\rho_{-}.$ Dans le cas du cercle, cette hypothèse n'est pas nécessaire.

L'hypothèse $\Gamma\Subset\Sigma\setminus\Phi$ permet de construire 
des quasimodes, puis à l'aide de ces derniers, de faire un problème de Grushin
pour ramener l'étude des valeurs propres à l'étude des zéros d'une fonction 
$\det E^{\delta}_{-+}$
(section 4 et 5). Après avoir rendu cette fonction holomorphe, nous concluons grâce à 
un lemme de comptage de zéros de fonctions holomorphes
(section 6).
La condition \ref{i3} sert à établir que $\det E^{\delta}_{-+}$ n'est pas trop petit 
avec une forte
probabilité.

Dans le cas classique, nous nous ramenons 
via une réduction semiclassique au cas précédent
(section 9) 
pour conclure avec le lemme de Borel-Cantelli afin d'avoir la loi de Weyl 
\textit{presque sûre} (section 10).

\paragraph{Rappel et notations.}
Précisons au préalable quelques notations, qui nous servirons par la suite.
Soit $m$ une fonction sur $\R$ de type $\langle \xi\rangle^{\ell},$ où $\ell\in\R$
 et $\Omega$ un ouvert de $\R^{2}.$ 
On introduit la classe $S(\Omega,m;\C^{n\times n})$
des symboles matriciels sur $\Omega$ 
\begin{align}
 \{ A(x,\xi)\in &C^{\infty}(\Omega,\C^{n\times n})|\,
\forall i,j,\;\forall\alpha,\beta\in\N,\; \exists C>0 \mbox { t.q.} \nonumber\\
& |\partial^{\alpha}_{x}
 \partial^{\beta}_{\xi}A^{i,j}(x,\xi)|\le C m(\xi) ,\; 
(x,\xi)\in \Omega\}.
\label{op.2}
\end{align}

Pour des symboles $A(x,\xi;h)$ dépendants de $h$, nous disons que $A\in S(m)$ 
si $A(.;h)$ est uniformément bornée dans $S(m)$ quand $h\in (0,1].$

Pour $k\in\R$, on pose $S^{k}(\Omega,m)=h^{-k}S(\Omega,m)$ et 
$S^{-\infty}(\Omega)=\bigcap S^{k}(\Omega,m).$

Soient $A,A_{k}\in S(\Omega,m),$  $k\ge 0.$ 
Si $\forall N\in\N,$  $A(x;h)-\sum_{0\le k\le N} A_{k}(x;h)h^{k}\in S^{-(N+1)}(\Omega,m),$
nous écrirons alors $A\sim \sum_{k=0}^{\infty}A_{k}h^{k}.$

Si $A$ et $B$ ont la même représentation asymptotique alors $A-B\in S^{-\infty}(\Omega,m).$

Si $A_{k}\in S(m),k\ge0$ alors il existe $A\in S(m)$ tel que 
$A\sim\sum A_{k}h^{k}.$

Un symbole $A\in S(m)$ est dit classique si $A\sim\sum A_{k}h^{k}$,
les fonctions matricielles $A_{k}$ étant indépendantes de $h$. $A_{0}$ est
dénommé le symbole principal de $A$. La classe des symboles classiques
est notée $S_{cl}(\Omega,m).$
\begin{prop}\label{op2}
L'application bilinéaire
\begin{eqnarray*}
S(\R^{2},m_{1})\times S(\R^{2},m_{2}) &\to & S(\R^{2},m_{1}m_{2})\\
(A_{1},A_{2})& \mapsto & A_{1}\# A_{2}
\end{eqnarray*}
où
\begin{equation}\label{op.4}
A_{1}\# A_{2}=e^{\frac{ih}{2}\sigma(D_{x},D_{\xi};D_{y},D_{\eta})}
A_{1}(x,\xi;h)A_{2}(y,\eta;h)|_{y=x, \eta=\xi}
\end{equation}
est continue. De plus, nous avons  la représentation asymptotique 
\begin{equation}\label{op.5}
(A_{1}\# A_{2})(x,\xi;h)\sim \sum_{k\ge 0}\frac{1}{k!}
\left(\frac{ih}{2}\sigma(D_{x},D_{\xi};D_{y},D_{\eta})\right)^{k}
A_{1}(x,\xi)A_{2}(y,\eta)\big|_{y=x,\eta=\xi}.
\end{equation}
\end{prop}

Grâce à (\ref{op.5}), il est possible de définir une composition pour les symboles définis
sur $\Omega,$
$S(\Omega,m_{1})\times S(\Omega,m_{2}) \to 
 S(\Omega,m_{1}m_{2})/S^{-\infty}(\Omega,m_{1}m_{2}).$
\begin{prop} \label{opd3}
Soit $A(x,\xi;h)\in S_{cl}(\Omega,m),$ les trois conditions suivantes sont équivalentes,

i) $A_{0}$ est inversible pour chaque $(x,\xi)\in\Omega$ et vérifie 
$A_{0}^{-1}=\mathcal{O}(\frac{1}{m}).$

ii) $A_{0}$ est inversible pour chaque $(x,\xi)\in \Omega$ et vérifie $A_{0}^{-1}\in S(m^{-1}).$

iii) $\exists B\in S(m^{-1})$ tel que 
\begin{eqnarray*}
A\# B&\sim&1\mbox{ dans }S(\Omega,1)\\
B\# A&\sim&1\mbox{ dans }S(\Omega,1).
\end{eqnarray*}
Un symbole qui vérifie i) est dit elliptique (au sens semiclassique).
\end{prop}

Lorsque $\Omega=\R^{2}$, on associe à
$A\in S(m)$ un opérateur pseudodifféren\-tiel
$A^{w}$ continue de $\mathcal{S}^{n}\to\mathcal{S}^{n}$ et 
de $(\mathcal{S'})^{n}\to (\mathcal{S'})^{n},$
défini par 
\begin{equation}\label{op.6}
A^{w}u(x):=\frac{1}{2\pi h}\iint e^{\frac{i}{h}(x-y)\xi} A(\frac{x+y}{2},\xi)u(y)\;dyd\xi.
\end{equation}
Puisque nous avons $A^{w}=\left((A^{i,j})^{w}\right)_{1\le i,j\le n},$
si  $A_{i}\in S(\R^{2},m_{i})$  nous obtenons alors la formule de composition 
$A_{1}^{w}A_{2}^{w}=(A_{1}\# A_{2})^{w}:\,\mathcal{S}^{n}\to\mathcal{S}^{n},\,
(\mathcal{S'})^{n}\to (\mathcal{S'})^{n}.$
\begin{theo} Si $A\in S(\R^{2},1),$ alors 
$A^{w}: L^{2}(\R,\C^{n})\to L^{2}(\R,\C^{n})$
est bornée, et sa norme est majorée par une constante indépendante de $h$.
\end{theo}
\begin{lemme}
Soit $A\in S_{cl}(\R^{2},m), A^{i,j}\sim\sum_{k\ge 0}h^{k}A_{k}^{i,j},$
introduisons
\[\mathrm{Supp}(A^{i,j}):=\overline{\bigcup_{k}\mathrm{supp}(A^{i,j}_{k}}).\]
Prenons $\chi\in C^{\infty}_{0}(\R^{2}),$ indépendant de  $h,$  alors 
\[ \forall i,j,\quad \mathrm{Supp}(A^{i,j})\cap\mathrm{supp}(\chi)=\emptyset
\Rightarrow \|(A\#\chi)^{w}\|_{L^{2}(\R)\to L^{2}(\R)}=\mathcal{O}(h^{\infty}).\]
\end{lemme}

Dans le cas scalaire, on pourra aussi consulter \cite{DiSj}, \cite{EvZw},
et dans le cas matriciel \cite{BoGl}, \cite{De}.\\

\section{Résultats préliminaires et Quasimodes}
On se place dans le cadre semiclassique et on étudie l'opérateur différentiel elliptique
(Hypothèse \ref{i1}) non-autoadjoint dans $L^{2}(S^{1},\C^{n})$
défini dans l'introduction. Rappelons que $q_{z}(x,\xi)=\det(p(x,\xi)-z).$

Si $z_{0}$ une valeur propre simple de $p(x_{0},\xi_{0}),$ où 
$(x_{0},\xi_{0})\in T^{\ast}S^{1},$
alors il existe un voisinage $U\subset  T^{\ast}S^{1}$ de $(x_{0},\xi_{0})$ 
et une fonction $C^{\infty},$ 
$\lambda :U\to \C,$ tel que $\lambda(x,\xi)$ soit une  valeur propre simple de 
$p(x,\xi)$ pour tout $(x,\xi)\in U,$ vérifiant au point $(x_{0},\xi_{0})$,
$\lambda (x_{0},\xi_{0})=z_{0}.$ 
\begin{prop}\label{sp2}
Soit $z_{0}$ une valeur propre simple de $p(x_{0},\xi_{0}),$ où $(x_{0},\xi_{0})\in T^{\ast}S^{1},$
alors nous avons l'équivalence
\[\frac{1}{2i}\{q_{z_{0}}(.),\overline{q_{z_{0}}(.)}\}(x_{0},\xi_{0})>0
\iff \frac{1}{2i}\{\lambda,\bar{\lambda}\}(x_{0},\xi_{0})>0.\]
\end{prop}
\textbf{Preuve.}
$q_{z}(x,\xi)$ se met sous la forme $g(x,\xi,z)(z-\lambda(x,\xi)),$
où $g(x,\xi,z)$ est polynomiale en $z$ et ne s'annule pas au point $(x_{0},\xi_{0},z_{0})$. 
Il faut ensuite remarquer que si $a(x,\xi)=b(x,\xi)c(x,\xi)$ et vérifie au point
$\rho_{0}=(x_{0},\xi_{0}),\,  
a(\rho_{0})=c(\rho_{0})=0$
et $b(\rho_{0})\ne0,$ alors 
\[\frac{1}{2i}\{a,\bar{a}\}(\rho_{0})=
|b(\rho_{0})|^{2}
\frac{1}{2i}\{c,\bar{c}\}(\rho_{0}).\]
\hfill $\square$ \medskip
\begin{prop}\label{sp3}
Soient $\rho_{0}=(x_{0},\xi_{0})$, $z_{0}\in\sigma(p(\rho_{0}))$.\\
(a) Si $\dim\mathcal{N}(p(\rho_{0})-z_{0})\ge 2$ alors 
\[\frac{1}{2i}\{q_{z_{0}}(.),\overline{q_{z_{0}}(.)}\}(\rho_{0})=0.\]
(b) Si $\dim\mathcal{N}(p(z_{0})-z_{0})=1$ alors il existe des matrices
$r_{0},s_{0}$ inversibles telles que $r_{0}^{-1}(p(\rho_{0})-z_{0})s_{0}$ 
admet 0 comme valeur propre simple.
\end{prop}
\textbf{Preuve.} 
(a) Pour une base convenable de $\C^{n},$ les deux première colonnes
de la matrice $p(\rho_{0})-z_{0}$ s'annulent. On voit donc que 
$\det (p(\rho)-z_{0})=\mathcal{O}(|\rho-\rho_{0}|^{2}).$ \\
(b) Soit $e_{1},\ldots,e_{n}$ une base telle que $(p(\rho_{0})-z_{0})e_{1}=0.$
Soit 
\begin{equation}
\begin{array}{l}
f_{2}=(p(\rho_{0})-z_{0})e_{2}\\
\vdots\\
f_{n}=(p(\rho_{0})-z_{0})e_{n}
\end{array}
\end{equation}
et $f_{1}$  tel que $f_{1},\ldots,f_{n}$ soit une base. Alors pour les bases $e_{1},\ldots,e_{n},$ et 
$f_{1},\ldots,f_{n}$
la matrice de $p(\rho_{0})-z_{0}$ devient
\begin{equation}\label{sp.5}
\left( \begin{array}{cccc}
0&0&\ldots&0\\
0&1& & \\
\vdots& &\ddots &\\
0&  & &1
\end{array}
\right).
\end{equation}
Il existe donc deux matrices de passage $r_{0},s_{0}$ pour lesquelles 
$r_{0}^{-1}(p(\rho_{0})-z_{0})s_{0}$ s'écrit comme dans (\ref{sp.5}).
\hfill $\square$ \medskip

Soit $z$ donné dans $\Sigma\setminus\Phi,$ l'image réciproque de $0$ par 
$q_{z}$ est donné par 
\begin{equation}
q^{-1}_{z}(0)=\{\rho_{+}^{\nu}(z), \rho_{-}^{\kappa}(z)|\,
\nu=1,\dots,\beta(z),
\kappa=1,\dots,\gamma(z)\}
\end{equation}
avec
\begin{equation}
\pm \frac{1}{2i}\{q_{z}(.),\overline{q_{z}(.)}\}(\rho_{\pm})>0.
\end{equation}
\begin{prop}\label{sp4}
a) Pour chaque $ z\in \Sigma\setminus\Phi,$ nous avons
$\beta(z),\gamma(z)<+\infty.$\\
b) Pour tout $ z\in \Sigma\setminus\Phi,$ nous avons $\beta(z)=\gamma(z).$\\
c) Si $z_{1},z_{2}$ appartiennent à une même composante  
connexe de $\Sigma\setminus\Phi,$ alors $\beta(z_{1})=\beta(z_{2}).$
\end{prop}
\textbf{Preuve.} Pour $a)$ et $c)$ c'est clair. 
$z_{0}$ étant fixé, on prend $q(x,\xi)\equiv q_{z_{0}}(x,\xi).$
On suppose pour se fixer les idées que $q(0,\xi)\ne 0,$
pour tout $\xi\in \R$ : il n'y a donc pas de points $\rho_{+}$ ou $\rho_{-}$
au dessus de $0.$ 
Nous coupons  le cercle $S^{1}\simeq \R/2\pi \Z,$ pour identifier,
avec l'application $(x,\xi)\mapsto (\mathrm{Re}\,w,\mathrm{Im} \,w),$ le tube 
$(S^{1}\setminus\{0\}) \times \{|\xi|\le C\}$ à un  rectangle $K$ de
$\C.$  Concrètement, nous pouvons  écrire
\begin{align*}
\partial K=
\underbrace{\{\xi=-C,x\in [0,2\pi]\}}_{\gamma_{1}}&\cup 
\underbrace{\{x=2\pi,|\xi |\le C\}}_{\gamma_{2}}\cup
\underbrace{\{\xi=C,x\in [0,2\pi]\}}_{\gamma_{3}}\\
&\cup \underbrace{\{x=0,|\xi |\le C\}}_{\gamma_{4}}.
\end{align*}

Puis nous  calculons 
la variation de l'argument de $q$ le long de la frontière de $K$
dans le sens positif.
Premièrement, puisque  
$q(x,\xi)=a(x)\xi^{mn}+\mathcal{O}(\xi^{mn-1}),$
avec 
$a(x)=\det(A_{m,0}(x))$ pour $\xi$ grand,
nous voyons que pour $C$ assez grand
\begin{align*}
\mathrm{var}\arg_{\gamma_{1}}\,q&=\mathrm{var}\arg_{S^{1}}a(x)\\
&=-\mathrm{var}\arg_{\gamma_{3}}\,q.
\end{align*}
Deuxièmement, comme $q(x,\xi)=q(x+2\pi,\xi),$ nous avons 
\[\mathrm{var}\arg_{\gamma_{2}}q+\mathrm{var}\arg_{\gamma_{4}}q=0.\]

Nous avons donc montré que la variation de l'argument de $q$ le long de $\partial K$ 
est nulle. Après une déformation de contour, nous pouvons aussi  écrire,  pour
$\epsilon$ assez petit, que
\begin{equation}
\mbox{var}\arg_{\partial K}q(x,\xi)=\sum_{\zeta\in q^{-1}(0)}
\mbox{var}\arg_{\partial D(\zeta,\epsilon)}q(x,\xi).
\end{equation}
On conclut alors avec le lemme qui suit :
\begin{lemme} 
Soit $q(\zeta)$ une fonction sur $\C_{\zeta}\equiv \R_{x}+i\R_{\xi},$ et $C^{\infty}$ 
dans un voisinage de $0.$
Si 
\begin{equation}\label{sp.9}
q(0)=0,\quad
\pm\frac{1}{2i}\{q,\bar{q}\}(0):=\pm\frac{1}{2i}(
\partial_{\xi}q\partial_{x}\bar{q}-
\partial_{x}q\partial_{\xi}\bar{q})(0)>0,
\end{equation} 
alors pour $\epsilon$
assez petit 
$\mathrm{var}\arg_{\partial D(0,\epsilon)} q(\zeta)=\pm 2\pi.$
\end{lemme}
\textbf{Preuve.} 
On fait un développement de Taylor de $q$ au voisinage  de zéro
\[q=a\,(\xi+ix)+b\,(\xi-ix)+\mathcal{O}(\|(x,\xi)\|^{2}),
\quad a,b\in\C. \]
Nous obtenons alors
\begin{align*}
\frac{1}{2i}\{q,\bar{q}\}(0)&=
|a|^{2}\,\frac{1}{2i}\{\xi+ix,\xi-ix\}(0)
+|b|^{2}\,\frac{1}{2i}\{\xi-ix,\xi+ix\}(0)\\
&=|b|^{2}-|a|^{2}.
\end{align*}

Deux cas se présentent.
Si $\frac{1}{2i}\{q,\bar{q}\}(0)>0$ alors $|b|>|a|$ et on voit que
$\mbox{var}\arg\,q=\mbox{var}\arg\,(\xi-ix)=-2\pi.$ 
Si $\frac{1}{2i}\{q,\bar{q}\}(0)<0$ alors $|b|<|a|$ et on a
$\mbox{var}\arg\,q=\mbox{var}\arg\,(\xi+ix)=+2\pi.$ 
\hfill $\square\square$ \medskip

Nous donnons maintenant un résultat d'existence de quasimode pour un système
différentiel matriciel qui généralise celui établi dans le cas scalaire par M.~Zworski
\cite{Zw1}. 
\begin{prop}
Pour tout $z$ dans $\Sigma(p)$ et 
$\rho_{0}:=(x_{0},\xi_{0})$ dans $T^{\ast}S^{1}$ avec
\[z\in \sigma(p(\rho_{0})),\quad \frac{1}{2i}\{q_{z}(.),\overline{q_{z}(.)}\}(\rho_{0})>0,\]
il existe $0\ne u(h)\in L^{2}(S^{1})$ tel que 
\[\|(P(h)-z)u(h)\|=\mathcal{O}(h^{\infty})\|u(h)\|.\]
De plus $u(h)$ a la forme $\chi(x)a(x;h)e^{i\varphi(x)/h}$ où $\chi$ est une troncature 
à support dans un voisinage de $x_{0}.$ Aussi $a(x_{0};h)\ne 0$ et 
$\im \varphi(x_{0})=\xi_{0}.$ 
\end{prop}
\textbf{Preuve.}
On suppose pour commencer que $z=0$ est une valeur propre simple de
 $p(\rho_{0})$. On cherche à construire des solutions BKW,
 $e^{i\varphi (x)/h}a(x;h),$ satisfaisant
 \begin{align}
 &e^{-i\varphi/h}P e^{i\varphi/h}a=\mathcal{O}(h^{\infty}),\label{sp.11}\\
&a(x;h)\sim a_{0}(x)+ha_{1}(x)+\ldots \mbox{ dans }
C^{\infty}(S^{1},\C^{n}),\,a_{0}\ne 0.\nonumber
 \end{align}
 La phase doit remplir l'équation eikonale
$ \det p(x,\varphi'(x))=0.$
Soit $\lambda$ la valeur propre
simple $C^{\infty}$ de $p(x,\xi)$ définie dans un voisinage de $\rho_{0}$ et vérifiant
$\lambda(x_{0},\xi_{0})=0.$ 
$\lambda$ est analytique en $\xi,$ et puisque $\partial_{\xi}\lambda(\rho_{0})\ne 0,$
on peut trouver une fonction $\varphi'$ définie dans un voisinage de $x_{0},$
telle que $\lambda(x,\varphi'(x))=0=\det p(x,\varphi'(x)),$ et $\varphi'(x_{0})=\xi_{0}.$

Les termes $a_{n}$ satisfont, eux, un système récurrent d'équations (les équations 
de transport).
La procédure pour donner les expressions explicites des $a_{n}$
 est décrite dans \cite{Fe} (p.54)   pour l'opérateur
 $hD_{x}+A(x).$ Par ailleurs, pour le cas scalaire on pourra consulter 
 \cite{Ha2}.  
 
On  prendra  ensuite comme quasimode $\chi a e^{i\varphi/h},$ où $\chi$
est à support compact dans un voisinage de $x_{0}.$
Pour la normalisation, on procède comme dans le cas scalaire, 
en remarquant que $\{\lambda,\bar{\lambda}\}(\rho)/2i>0,$ ce qui entraine 
que $\mbox{Im } \varphi''(x_{0})>0,$ voir \cite{Da}.

Pour le cas où $z_{0}$ est une valeur propre multiple, on est ramené au cas d'une valeur
propre simple après composition par $r_{0}$ et $s_{0}$ (proposition \ref{sp3}).
\hfill $\square$ \medskip

Pour une étude plus approfondie de l'existence de quasimodes
pour les systèmes d'opérateurs semiclassiques, on consultera
\cite{De}.

\section{Problème de Grushin pour l'opérateur non-perturbé}
Pour le problème de Grushin, seule l'hypothèse d'ellipticité est imposée,
donc la condition \ref{i3} n'est pas nécessaire. Nous rappelons que
$q_{z}(x,\xi)$ désigne le déterminant de $p(x,\xi)-z.$ On se place dans le 
cadre semiclassique.

Soit $z_{0}$ un point de $\Lambda(p)=\Sigma\setminus\Phi,$ 
et prenons les points $\rho_{\pm}^{\nu}(z_{0})$
de $q^{-1}_{z_{0}}(0)$ tels que 
\begin{equation}\label{gg.1}
\pm\frac{1}{2i}\{q_{z_{0}}(.),\overline{q_{z_{0}}(.)}\}(\rho^{\nu}_{\pm}(z_{0}))>0.
\end{equation}

Nous omettons dans la suite d'écrire l'indice $\nu$.
Nous indiquerons dans le texte quand nous en aurons  besoin.

Comme $dq_{z_{0}}, d\bar{q}_{z_{0}}$ sont linéairement indépendant, 
il existe un voisinage $U(z_{0})$ de $z_{0}$ et 
$\rho_{\pm}(z)\in C^{\infty}(U(z_{0}))$
pour lesquel $q(\rho_{\pm}(z),z)=0,$ et 
\begin{equation}
\pm\frac{1}{2i}\{q_{z}(.),\overline{q_{z}(.)}\}(\rho_{\pm}(z))> 0.
\end{equation}

On suppose pour commencer que $0$ est une valeur propre simple de 
$p(\rho_{\pm}(z_{0}))-z_{0}.$
Dans le cas général $\mbox{dim }\mathcal{N}(p(\rho_{\pm}(z_{0}))-z_{0})=1,$ 
la proposition 
\ref{sp4} montre qu'après composition par $r_{0}$ et $s_{0}$ on 
est ramené au cas d'une valeur propre simple.

Utilisant le paragraphe 3 du Ch.I de \cite{GoKr}, 
on déduit qu'il existe un voisinage $W(z_{0})$ de $z_{0}$, un voisinage 
$V_{\pm}$ de $\rho_{\pm}(z_{0})$ pour lesquels il existe une matrice
$u_{\pm}(x,\xi)$ inversible pour chaque $(x,\xi)\in V_{\pm}$, une fonction 
scalaire  $\lambda_{\pm}(x,\xi)$ 
vérifiant 
\begin{equation}\label{gg.2}
\lambda_{\pm}(\rho_{\pm}(z))=z,\quad\forall(x,\xi)\ne\rho_{\pm}(z),\;
\lambda_{\pm}(x,\xi)\ne z,
\end{equation}
et une matrice $(n-1)\times (n-1)$, $h_{\pm}(x,\xi),$ avec 
\begin{equation}\label{gg.3}
\forall (x,\xi)\in V_{\pm},
\forall z\in W(z_{0}),\quad\det(h_{\pm}-z)\ne 0,
\end{equation}
tels que 
\begin{multline}\label{gg.4}
\forall(x,\xi)\in V_{\pm},\forall z\in W(z_{0}),\\
u_{\pm}(x,\xi)(p(x,\xi)-z)u^{-1}_{\pm}(x,\xi)=
\left(\begin{array}{ccc}
\lambda_{\pm}(x,\xi) -z&0 \\
 0& h_{\pm}(x,\xi )-z
\end{array}\right).
\end{multline}

Puisque $V_{\pm}$ est relativement compact,
$u_{\pm},u^{-1}_{\pm},h_{\pm}\in S(V_{\pm},1)$
et $p\in  S(V_{+},1),$ $S(V_{-},1)$.
On adap\-te ensuite un résultat du à M. Taylor \cite{Ta}, 
voir aussi \cite{HeSj}, proposition 3.1.1,
\begin{prop}\label{gg1}
Soit $\Omega$ un ouvert de $T^{\ast}S^{1}.$
Soit $A\in S_{cl}(\Omega,1)$, dont le symbole principal vérifie
\[U_{0}A_{0}U^{-1}_{0}=
\widetilde{A}_{0}:=
\left(\begin{array}{ccc}
\widetilde{A}_{0}^{1,1}&0 \\
 0& 
 \widetilde{A}_{0}^{2,2}
\end{array}\right), 
\quad U_{0},U^{-1}_{0}\in S(\Omega,1),\]
où pour chaque $(x,\xi)$,
$\widetilde{A}_{0}^{1,1}(x,\xi)$ et $\widetilde{A}_{0}^{2,2}(x,\xi)$ 
ont des spectres disjoints.
Il existe alors $U\in S_{cl}(\Omega,1),$
$\widetilde {U}\in S_{cl}(\Omega,1)$ vérifiant
 \[U\#\widetilde{U}\sim 1,\quad
\widetilde{U}\#U\sim 1,\quad
 U= U_{0}\mathrm{mod}\,h S_{cl}(\Omega,1),\]
dans $S(\Omega,1)$  tels que 
\[U\# A\#\widetilde{U}\sim
\left(\begin{array}{ccc}
\widetilde{A}^{1,1}& 0\\
0 &\widetilde{A}^{2,2}
\end{array}\right), \]
où le symbole principal, 
$\widetilde{A}^{i,i}_{0}$, de $\widetilde{A}^{i,i}$ vérifie $\widetilde{A}^{i,i}_{0}=A^{i,i}_{0},$
$i=1,2.$
\end{prop}

Dans notre cas, nous obtenons:
\begin{corol}\label{gg2}
Soit $P(x,\xi)-z$ le symbole de l'opérateur $P-z$, alors il existe  $U_{\pm},\widetilde{U}_{\pm},
H_{\pm}$, et  $\tilde{\lambda}_{\pm}\in S_{cl}(V_{\pm},1)$ tels que
\[ U_{\pm}(x,\xi;h)\# (P(x,\xi)-z)\# \widetilde{U}_{\pm}(x,\xi;h)\sim
\left(\begin{array}{ccc}
\tilde{\lambda}_{\pm}(x,\xi; h)-z&0 \\
0 & H_{\pm}(x,\xi; h)-z
\end{array}\right)\]
dans $S(V_{\pm},1),$ où le symbole pricipal de 
$\tilde{\lambda}_{\pm}$ est $\lambda_{\pm}$ et celui de $H_{\pm}$, 
$h_{\pm}.$
\end{corol}

$H_{\pm}$ est elliptique au sens semiclassique, et 
compte tenu de  la proposition \ref{sp2}, le symbole principal  $\lambda_{\pm}$ de 
$\tilde{\lambda}_{\pm}$
vérifie pour chaque $z\in W(z_{0}),$
\begin{equation}\label{gg.5}
\rho_{\pm}(z)\in \lambda_{\pm}^{-1}(z),\quad \pm\frac{1}{2i}
\{\lambda_{\pm},\bar{\lambda}_{\pm}\}(\rho_{\pm},z)>0.
\end{equation}
\textbf{On est ainsi ramené au cas scalaire traité dans \cite{Ha2}}.
La proposition 3.3 de \cite{Ha2} montre qu'il existe 
un voisinage $\widetilde{W}(z_{0})\subset W(z_{0})$ de $z_{0}$, un voisinage 
$\widetilde{V}_{\pm}\subset V_{\pm}$ 
contenant $\rho_{\pm}(z)=(x_{\pm}(z),\xi_{\pm}(z))$, $z\in\widetilde{W}(z_{0})$, deux symboles
$q_{\pm}\in S_{cl}(\widetilde{V}_{\pm},1),$ et  $g_{\pm}\in 
S_{cl}(\pi_{x}(\widetilde{V}_{\pm}),1)$
qui dépendent de manière $C^{\infty}$ de $z\in\widetilde{W}(z_{0})$ tels que
\begin{align}
\tilde{\lambda}_{+}(x,\xi;h)-z&\sim q_{+}(x,\xi,z;h)\#
(\xi+g_{+}(x,z;h)) \mbox{ dans } S(\widetilde{V}_{+},1),\label{gg.6}\\
\tilde{\lambda}_{-}(x,\xi;h)-z&\sim (\xi+g_{-}(x,z;h))\#
q_{-}(x,\xi,z;h) \mbox{ dans } S(\widetilde{V}_{-},1),\label{gg.7}
\end{align}
avec $q_{\pm, 0}(\rho_{\pm}(z),z)\ne 0,$ et  
ou $g_{\pm,0}(x_{\pm}(z),z)=-\xi_{\pm}(z)$ pour
$z\in\widetilde{W}(z_{0}).$

La fonction $g_{\pm}$ est définie sur $\pi_{x}(\widetilde{V}_{\pm})$. On prolonge 
$g_{\pm}$ dans $C^{\infty}(\R)$ de telle sorte que 
\begin{equation}
g_{\pm}(y)=\mp \frac{i}{C_{\pm}}(y-x_{\pm}),\quad
|y|\ge C,\quad C_{\pm}>0,
\end{equation}
et aussi $\mbox{Im }g_{\pm}(y)\ne 0$ pour $y\ne x_{\pm}(z).$

\textbf{On identifiera fréquemment les intervalles de $\R$ de longueur
$<2\pi$ à des intervalles de $S^{1}.$}

Soit $\Upsilon_{\pm}\in L^{2}(\R)$ les solutions normalisées de
\begin{align*}
(hD_{x}+g_{+})\widetilde{\Upsilon}_{+}&=0,\\
 (hD_{x}+g_{-})^{\ast}\widetilde{\Upsilon}_{-}&=0.
 \end{align*}
$\widetilde{\Upsilon}_{\pm}$ est de la forme (type BKW)
$h^{-\frac{1}{4}}a_{\pm}(x,z;h)e^{\frac{i}{h}\varphi_{\pm}(x, z)},$
où $a_{+}$ est un symbole classique avec $a_{\pm,0}(x_{\pm}(z),z)=0$, et
\begin{equation}
\varphi'_{\pm}(x_{\pm}(z),z)=0, \;
\partial_{x}\varphi_{\pm}(x,z)=g_{\pm,0}(x,z),\,
\im\varphi_{\pm}\ge0.
\end{equation}
Puisque $\partial_{x}\varphi_{\pm}(x,z)=g_{\pm,0}(x,z)$ et 
$\lambda_{\pm}(x,\xi)-z=q_{\pm,0}(\xi+g_{\pm,0})$ alors
\begin{equation}
\lambda_{\pm}(x,-\partial_{x}\varphi_{\pm}(x,z))-z=0,
\quad \forall x\in \widetilde{V}_{\pm},\;z\in \widetilde{W}(z_{0}).
\end{equation}
On note $\Upsilon_{\pm}:=(\widetilde{\Upsilon}_{\pm},0,\ldots,0)\in C^{\infty}(\R;\C^{n}).$

On réintègre l'indice $\nu.$

Sans perte de généralité on suppose que $z_{0}$ est une valeur propre simple
de $p(\rho^{\nu}_{\pm}(z_{0}))$ pour tout $1\le\nu\le \beta,$
$\beta$ étant la valeur constante de $\beta(z)$ sur une composante connexe 
de $\Sigma\setminus\Phi.$ 

Il existe un voisinage de $z_{0}$ noté
$\widetilde{W}(z_{0})\subset\widetilde{W}^{\nu}(z_{0})$ pour tout $\nu,$
et 
des voisinages $\Theta_{\pm}^{\nu}\subset \widetilde{V}^{\nu}_{\pm} $ 
contenant $\rho^{\nu}_{\pm}(z)$ pour $z\in \widetilde{W}(z_{0}),$ et tels que\
$\Theta_{\pm}^{\nu}$ sont à adhérences disjointes, ce qui signifie  
que pour $i,j=\pm$, 
$\overline{\Theta^{\nu}_{i}}\cap\overline{\Theta^{\kappa}_{j}}
=\emptyset$ sauf pour 
$(i=j=\pm$ avec $\nu=\kappa).$ 

Définissons alors les fonctions suivantes
\begin{equation}
\chi_{\pm}^{\nu}\in C^{\infty}_{0}(\Theta_{\pm}^{\nu}) \mbox{ indépendants de } h,
\; \chi_{\pm}^{\nu}=1\mbox{ près de } \rho^{\nu}_{\pm}(\widetilde{W}(z_{0})),
\label{gg.8}\\
\end{equation}
\begin{equation*}
\phi^{\nu}_{\pm}\in C^{\infty}_{0}(\pi_{x}(\Theta^{\nu}_{\pm})),\mbox{ indépendants de } h,
\;\phi_{\pm}^{\nu}=1\mbox{ près de } \pi_{x}(\rho^{\nu}_{\pm}(\widetilde{W}(z_{0})).
 \label{g.}
\end{equation*} 
Soient $\widehat{\chi}_{\pm}^{\nu}
\in C^{\infty}_{0}(\Theta^{\nu}_{\pm})$ et 
$\widetilde{\phi}_{\pm}^{\nu}
\in C^{\infty}_{0}(\pi_{x}(\Theta^{\nu}_{\pm}))$
satisfaisant
$\chi_{\pm}^{\nu}\prec\widehat{\chi}_{\pm}^{\nu}$
et
$\phi_{\pm}^{\nu}\prec\widetilde{\phi}_{\pm}^{\nu}.$
Introduisons également  les fonctions suivantes, définies sur $S^{1},$
\begin{align}
e_{+}^{\nu}&=\phi_{+}^{\nu}(\widetilde{U}_{+}^{\nu}\#\chi_{+}^{\nu})^{w} \Upsilon_{+}^{\nu},\quad
&f_{+}^{\nu}&=\widetilde{\phi}_{+}^{\nu}((U^{\nu}_{+})^{\ast}\#\widehat{\chi}_{+}^{\nu})^{w}\Upsilon_{+}^{\nu},\\
e_{-}^{\nu}&=\phi_{-}^{\nu}(\chi_{-}^{\nu}\#(U^{\nu}_{-})^{\ast})^{w}\Upsilon_{-}^{\nu},\quad
& f_{-}^{\nu}&=\widetilde{\phi}_{-}^{\nu}(\widehat{\chi}_{-}^{\nu}\#\widetilde{U}_{-}^{\nu})^{w}\Upsilon_{-}^{\nu}.
\end{align}

\textit{Rappelons que les multiplications par les troncatures  
$\phi$ et $\widetilde{\phi}$ servent à identifier les intervalles de longueur $<2\pi$ de $\R$
(où agissent nos opérateurs pseudodifférentiels) et de $S^{1}$. }

Etant donné que $\widetilde{\Upsilon}_{\pm}^{\nu}$ est une fonction de type BKW
microlocalisée près de $\rho_{\pm}(z),$ nous avons 
\begin{equation}\label{gg.11}
\langle e_{\pm}^{\nu},f_{\pm}^{\nu}\rangle=1+\mathcal{O}(h^{\infty}), 
\mbox{ pour tout } \nu, 
\end{equation}
de plus il existe une constante $C>0,$ indépendante de $h$,
telle que $\frac{1}{C}\le\| e_{\pm}^{\nu}\|,\|f_{\pm}^{\nu}\|\le C.$
Normalisons $e_{\pm}^{\nu}$ et multiplions en conséquence 
$f_{\pm}^{\nu}$ 
par une constan\-te minorée et majorée uniformément
par rapport à $h,$ pour que (\ref{gg.11}) reste vérifié.
\begin{prop}\label{gg5}
$e_{\pm}^{\nu}\in L^{2}(S^{1})$ sont des fonctions normalisées 
de type BKW vérifiant
\begin{equation}\label{gg.12}
\|(P-z)e_{+}^{\nu}\|_{L^{2}(S^{1})},
\|(P-z)^{\ast}e_{-}^{\nu}\|_{L^{2}(S^{1})}=\mathcal{O}(h^{\infty}).
\end{equation}
De plus, 
$e_{\pm}^{\nu}$ admet une expression de la forme
\begin{equation}\label{gg.13}
e^{\frac{i}{h}\varphi_{\pm}^{\nu}(x,z)}I_{\pm}^{\nu}(x,z;h)+
r_{\pm}^{\nu}(x;h)
,\quad  I_{\pm}\in C^{\infty}(S^{1},\C^{n}),
\end{equation}
où les coefficients $I_{\pm}^{\nu,k}(x,z;h)$ du vecteur $I_{\pm}^{\nu}(x,z;h)$
ne s'annulent pas au point $x_{\pm}^{\nu}(z),$ et admettent un
développement en puissances de $h,$ dans $C^{\infty}(S^{1}),$
de la forme
\begin{equation}\label{gg.14}
I^{\nu,k}_{\pm}(x,z;h)\sim  h^{-1/4}(I^{\nu,k}_{\pm,0}(x,z)+hI^{\nu,k}_{\pm,1}(x,z)+\ldots)
\end{equation}
et où $\varphi_{\pm}$ a été introduit pour $\widetilde{\Upsilon}_{\pm},$
et $r^{\nu}(x;h)\in S^{-\infty}(1),$ pour tout $\nu.$ 
$e_{\pm}^{\nu}$ est donc microlocalement concentré
près de  $\rho_{\pm}^{\nu}(z)=(x_{\pm}^{\nu}(z),\xi_{\pm}^{\nu}(z)).$ 
$f_{\pm}^{\nu}$ admet une représentation similaire
à $e_{\pm}^{\nu}.$
\end{prop}
L'expression (\ref{gg.13}) résulte d'un résultat de Melin-Sj\"ostrand 
\cite{MeSj}
sur l'action d'un opérateur pseudodifférentiel sur une fonction 
BKW avec une phase complexe $\varphi$ admettant un point critique 
non-dégénéré et vérifiant $\im \varphi\ge 0.$
\begin{theo}\label{gg6} 
Pour tout $z$ dans $\widetilde{W}(z_{0}),$
\[ \mathcal{P}=
\left(
\begin{array}{ccc}
P-z  & R_{-}\\
R_{+} & 0
\end{array}
\right): H^{m}_{sc}(S^{1};\C^{n})\times \C^{\beta}\to 
L^{2}(S^{1};\C^{n})\times \C^{\beta},
\]
avec 
\begin{align}
(R_{+}u)_{\nu}:=\langle u,f_{+}^{\nu}\rangle,\;u\in H^{m}_{sc},\\
 R_{-}u_{-}:=\sum_{1\le\nu\le \beta}
f^{\nu}_{-}u^{\nu}_{-},\; u_{-}\in \C^{\beta},
\end{align}
est inversible d'inverse
\[ \mathcal{E}=
\left(
\begin{array}{ccc}
E& E_{+}\\
E_{-}& E_{-+}
\end{array}
\right)
=
\left(
\begin{array}{ccc}
E_{0}+\mathcal{O}(h^{\infty})& F_{+}+\mathcal{O}(h^{\infty})\\
G_{-}+\mathcal{O}(h^{\infty})& \mathcal{O}(h^{\infty})
\end{array}
\right)
\]
où $E_{0}=\mathcal{O}(\frac{1}{\sqrt{h}}),$ et
\begin{align}
F_{+}v_{+}=\sum_{1\le\nu\le \beta} v_{+}^{\nu}e_{+}^{\nu},\;v_{+}\in \C^{\beta},\\
(G_{-}v)_{\nu}=\langle v,e^{\nu}_{-}\rangle, \;
v\in L^{2}.
\end{align}

De plus $E_{0}$ ne propage pas les supports au sens où si 
$\psi_{1},\psi_{2}\in S(T^{\ast}\R,1)$ sont à  support disjoint
avec $|\pi_{x}(\mathrm{supp}\,\psi_{k})|<2\pi,$ alors pour tous
$\chi_{k}(x)\prec \widetilde{\chi}_{k}(x)\in C^{\infty}_{0}(S^{1}),$
$k=1,2$ avec $\{x;\,\pi_{x}(\psi_{k})=1\}\subset \mathrm{supp}\,\chi_{k}<2\pi,$
alors 
\[(\widetilde{\chi}_{2}\psi_{2}^{w}\chi_{2})
E_{0}
(\widetilde{\chi}_{1}\psi_{1}^{w}\chi_{1})
=\mathcal{O}(h^{\infty}),\mbox{ dans }
\mathcal{L}(L^{2}(S^{1}),H^{1}_{sc}(S^{1})).\]
Cette dernière propriété est à relier aux lemmes 4.6 et 4.8 de Hager \cite{Ha2}.
\end{theo}
\section{Problème de Grushin pour l'opérateur perturbé}
 Commen\c cons par rappeler la  proposition suivante qui améliore
 celle de \cite{HaSj}, section 6. Pour la preuve qui simplifie celle de 
 \cite{HaSj}, on consultera \cite{Bo}. 
\begin{prop}\label{pp1}
Soit $(Y_{k})_{k\in\Z}$ une suite de variables complexes indépendantes
de loi $Y_{k}\sim \mathcal{N}(0,\sigma^{2}_{k}).$ Si 
$\sum\sigma_{k}^{2}<\infty,$ alors on a
\[\forall x>0,\quad \mathbb{P}(\sum_{k\in\Z}|Y_{k}|^{2}\ge x)
\le \exp[\frac{C_{0}}{2s_{1}}\sum_{k\in\Z}\sigma^{2}_{k}-\frac{x}{2s_{1}}].\]
Ici $s_{1}=\max \sigma_{k}^{2},$ et $C_{0}>0.$
\end{prop}
\begin{corol}\label{pp2}
Soit $(Y_{k})_{k\in\Z}$ une suite de variables complexes indépendantes
de loi $Y_{k}\sim \mathcal{N}(0,\sigma^{2}_{k}).$ Si 
$\sum\sigma_{k}<\infty,$ alors il existe $C_{0}>0$ tel que 
\[\forall x>0,\quad \mathbb{P}(\sum_{k\in\Z}|Y_{k}|\ge x)
\le \exp\left[\frac{C_{0}}{2\|\sigma\|_{\infty}}\|\sigma\|_{1}-\frac{x^{2}}{2\|\sigma\|_{\infty}\|\sigma\|_{1}}
\right],\] 
où $\|\cdot\|_{p}$ désigne la norme $\ell^{p}.$ 
\end{corol}
\textbf{Preuve.} 
Par Cauchy-Schwarz
\[\mathbb{P}\left(\sum_{k\in\Z}|Y_{k}|\ge x\right)\le
\mathbb{P}\left( (\sum_{k\in \Z} \sigma_{k})^{\frac{1}{2}}
\left(\sum_{k\in\Z} |\sqrt{\sigma_{k}} X_{k}|^{2}\right)^{1/2}
\ge x\right),\]
où $Y_{k}=\sqrt{\sigma_{k}}X_{k},$ $X_{k}\sim\mathcal{N}(0,1).$ On utilise la proposition \ref{pp1} pour achever
la preuve. 
\hfill $\square$ \medskip

La norme de $\|Q_{\omega}\|_{H_{sc}^{m}(S^{1})\to L^{2}(S^{1})}$ est majorée par 
\begin{equation}
C\sup_{\alpha_{0}\le\alpha\le\alpha_{1},\, x\in S^{1}}\|Q_{\alpha}(x)\|\le
\widetilde{C}\sup_{\alpha,i,j} \|Q^{i,j}_{\alpha}(x)\|_{L^{\infty}},
\end{equation}
 où $C,\widetilde{C}$ sont des constantes strictement positives et 
$\|\cdot\|$ est une norme sur $\C^{n}\times\C^{n}.$

\begin{prop}\label{pp3}
On suppose admis l'hypothèse \ref{i4}.
Il existe $C>0$ tel que pour chaque $x>0,$ et $0<h\ll1,$ on ait
\begin{equation}\label{pp.4}
\mathbb{P}(\|Q_{\omega}\|_{H^{m}_{sc}(S^{1})\to L^{2}(S^{1})}\le x)
\ge 1-
\exp(C-\frac{x^{2}}{C}).
\end{equation}
$Q_{\omega}$ est 
donc bornée presque sûrement comme opérateur de
$H^{m}_{sc}(S^{1})\to L^{2}(S^{1}).$
\end{prop}
\noindent\textbf{Preuve.}
Majorant $\|Q^{i,j}_{\alpha}(x)\|_{L^{\infty}}$ par 
la somme des valeurs absolues $\sum_{k\in\Z} |q^{i,j}_{\alpha,k}|,$ nous avons 
\begin{align}\label{pp.4bis}
\mathbb{P}(\|Q\|_{H^{m}_{sc}(S^{1})\to L^{2}(S^{1})}\ge x)
&\le 
\mathbb{P}(\sum_{\alpha,i,j}\|Q^{i,j}_{\alpha}(x)\|_{L^{\infty}}\ge \frac{x}{\widetilde{C}})\\
&\le
\mathbb{P}(\sum_{\alpha,i,j,k}|q^{i,j}_{\alpha,k}|\ge \frac{x}{\widetilde{C}}).
\nonumber
\end{align}
Nous utilisons le corollaire \ref{pp2}, (\ref{pp.4bis}) devient alors
\[\le \exp [\frac{C_{0}}{2\sup \sigma^{i,j}_{\alpha,k}}\sum \sigma^{i,j}_{\alpha,k}
-\frac{x^{2}}{2\widetilde{C}^{2}\,\sup (\sigma^{i,j}_{\alpha,k})\sum \sigma^{i,j}_{\alpha,k}}],
\]
ce qui termine la preuve puisque par hypothèse $\sigma_{k}$ est sommable.
\hfill $\square$ \medskip

\begin{corol}\label{pp4} On suppose que l'hypothèse \ref{i4} est vérifiée.
Il existe alors $C>0$  tel que pour tout $0<h\ll 1$  
\[\mathbb{P}(\|Q_{\omega}\|_{H^{m}_{sc}(S^{1})\to L^{2}(S^{1})}\le \ln (h^{-1}))\ge
1-Ce^{-\frac{1}{C}(\ln h)^{2}}.\]
\end{corol}

Dans la suite on travaille sous l'hypothèse \ref{i4}.
\begin{prop}\label{pp5} 
Soit $\delta\ll \sqrt{h}|\ln h|^{-1}$ un paramètre de perturbation 
et $z_{0}$ dans $\Sigma\setminus\Phi.$
Il existe un voisinage $\widetilde{W}(z_{0})$ de $z_{0}$ inclus dans 
$\Sigma\setminus\Phi$ tel que
avec une probabilité 
\[\ge 1-Ce^{-\frac{1}{C}(\ln h)^{2}}\]
pour tout $z$ dans $\widetilde{W}(z_{0}),$
\[\mathcal{P}^{\delta}=
\left(\begin{array}{ccc}
P-z-\delta Q & R_{-}\\
R_{+} & 0 
\end{array}
\right)
\]
est continu $H^{m}_{sc}(S^{1})\times \C^{\beta }\to L^{2}(S^{1})\times \C^{\beta}$ 
et admet 
un inverse $\mathcal{E}^{\delta}$
de la forme
\begin{align}
\mathcal{E}^{\delta}
&=\mathcal{E}^{0}+
\left(
\begin{array}{ccc}
\sum_{j\ge 1}E(\delta QE)^{j} & \sum_{j\ge 1}(E\delta Q)^{j}E_{+}\\
\sum_{j\ge 1} E_{-}(\delta QE)^{j} & \sum_{j\ge 1} E_{-}(\delta QE)^{j-1}
(\delta QE_{+})\\
\end{array}
\right)\label{pp.5}\\
&=\mathcal{E}^{0}+
\ln (h^{-1})
\left(\begin{array}{ccc}
\mathcal{O}(\frac{\delta}{h})& 
\mathcal{O}(\frac{\delta}{\sqrt{h}})\\
\mathcal{O}(\frac{\delta}{\sqrt{h}})&
\mathcal{O}(\delta)\\
\end{array}
\right).\label{pp.6}
\end{align}
\end{prop}
\textbf{Preuve.}
Nous avons $\mathcal{P}^{\delta}\mathcal{E}=1-K$ où
\[K=\left(\begin{array}{ccc}
\delta QE & \delta QE_{+}\\
0 & 0
\end{array}
\right).
\]

Il existe $C>0$ tel que avec une probabilité supérieure à
$1-Ce^{-\frac{1}{C}(\ln h)^{2}}$ on ait 
$\|Q\|_{H^{m}_{sc}\to L^{2}}\le|\ln h|,$ impliquant 
\[\|K\|\lesssim\delta\|Q\|\|E\|=\delta h^{-\frac{1}{2}}|\ln |\ll 1.\]
\hfill $\square$ \medskip

Dans la suite, on suppose que $\|Q\|_{H^{m}_{sc}\to L^{2}}\le|\ln h|,$ 
et $\delta\ll \sqrt{h}|\ln h|^{-1}.$\\

\section{Propriétés d'holomorphie de $E_{-+}$}
Puisque $\partial_{\bar{z}}(\mathcal{P}\mathcal{E})=0,$ nous avons 
\begin{equation}\label{pp.11}
\partial_{\bar{z}} E_{-+}=-E_{-+}(\partial_{\bar{z}}R_{+})E_{+}-
E_{-}(\partial_{\bar{z}}R_{-})E_{-+}.
\end{equation}
Donc, grâce à la cyclicité de la trace, nous obtenons
\begin{align}\label{pp.12}
\partial_{\bar{z}}\det E_{-+}&=\mbox{tr}(( \partial_{\bar{z}}E_{-+})E^{-1}_{-+})\det E_{-+}
\nonumber\\
&=-\mbox{tr}((\partial_{\bar{z}}R_{+})E_{+}+E_{-}(\partial_{\bar{z}}R_{-}))
\det E_{-+}\nonumber\\
&=:-k^{0}(z)\det E_{-+}(z).
\end{align}
Dès lors, si on choisit une solution de l'équation
\begin{equation}\label{pp.13}
\frac{1}{h}\partial_{\bar{z}}l^{0}=k^{0},\quad
l^{0}(z)=\frac{h}{\pi}\int_{\widetilde{W}(z_{0})}\frac{k^{0}(z')}{z-z'}
d\mathrm{Re}z'd\mathrm{Im}z'
\end{equation}
dans un voisinage de $\widetilde{W}(z_{0})$, nous obtenons une fonction 
$e^{l^{0}/h}\det E_{-+}$ holomorphe
avec les mêmes zéros que $\det E_{-+}$ dans $\widetilde{W}(z_{0}).$
\begin{prop}\label{pp6}
$\Delta\mathrm{Re}\,l^{0}(z) ,$ défini sur $\widetilde{W}(z_{0}),$ 
est strictement sousharmonique et 
\begin{equation}
(\Delta\mathrm{Re}\,l^{0}(z)+\mathcal{O}(h))\,
d \mathrm{Re}\, z\wedge d\mathrm{Im}\, z=
\sum_{1\le \nu\le \beta}(d\xi_{-}^{\nu}\wedge dx_{-}^{\nu}
-d\xi_{+}^{\nu}\wedge dx_{+}^{\nu}).
\end{equation}
\end{prop}

Dans \cite{Ha2}, Hager utilise des arguments géométriques pour démontrer ce résultat.
Nous proposons ici une preuve directe.

\noindent
\textbf{Preuve.}
Rappelons $\Delta:=4\partial_{z}\partial_{\bar{z}}.$ Nous avons montré que
\begin{align}
\frac{1}{h}\partial_{\bar{z}}l_{0}=k_{0}(z)
&= \sum_{1\le\nu\le \beta}((\partial_{\bar{z}}R_{+})E_{+})_{\nu,\nu}+
(E_{-}(\partial_{\bar{z}}R_{-}))_{\nu,\nu}\\
&=\sum_{1\le\nu\le\beta}\langle e_{+}^{\nu},\partial_{z} f_{+}^{\nu}\rangle+
\langle\partial_{\bar{z}} f_{-}^{\nu},
e_{-}^{\nu}\rangle+\mathcal{O}(h^{\infty}).
\end{align}
Nous avons
\[ h\langle e_{+},\partial_{z} f_{+}\rangle =
\langle e_{+}, i(\partial_{z}\varphi_{+})(x,z)\,f_{+}\rangle
+\mathcal{O}(h).\]
Puisque $\langle e_{+},f_{+}\rangle =1+\mathcal{O}(h^{\infty}),$ alors le lemme de la phase
stationnaire implique que 
\[ h\langle e_{+},\partial_{z} f_{+}\rangle=
-i\,\overline{(\partial_{z}\varphi_{+})(x_{+}(z),z)}+\mathcal{O}(h).\]
De même, on montre que
\[h\langle \partial_{\bar{z}}f_{-},e_{-}\rangle=
i(\partial_{\bar{z}}\varphi_{-})(x_{-}(z),z)+\mathcal{O}(h).\]

Il ne reste alors plus qu'à achever la preuve avec le lemme suivant:
\begin{lemme}\label{pp7}
Soit $\varphi_{+}(x,z)$ dans $C^{\infty}(\R\times \C)$ et vérifiant pour tous $z$
\begin{align}
\varphi_{+}(x_{+}(z),z)&=0,\quad (\partial_{x}\varphi_{+})(x_{+}(z),z)=\xi_{+}(z)\in\R,
\end{align}
et dans un voisinage de $x_{+}(z)$
\begin{align}
\lambda_{+}(x_{+},\varphi'_{x}(x,z))-z=0,
\end{align}
alors pour tout $z,$ 
\begin{equation}\label{t.21}
4\,\mathrm{Im }\;\frac{\partial}{\partial z}\bigg(\overline{(\partial_{z}
\varphi_{+})(x_{+}(z),z)}\bigg)(z) \,d\mathrm{Re}\,z\wedge d\mathrm{Im}\,z= -
d\xi_{+}\wedge dx_{+}.
\end{equation}
\end{lemme}
On a un lemme similaire avec $\varphi_{-}$.\\
\textbf{Preuve.} 
Commen\c cons par remarquer que, avec $\lambda\equiv \lambda_{+}$ et 
$\varphi\equiv\varphi_{+}$: 
\begin{equation}
\left(
\begin{array}{ccc}
\lambda'_{x}& \lambda'_{\xi}\\
\bar{\lambda}'_{x}&\bar{\lambda}'_{\xi}
\end{array}\right)
\left(\begin{array}{c}
\partial_{z}x_{+}\\
\partial_{z}\xi_{+}
\end{array}\right)
=
\left(\begin{array}{c}
1\\
0
\end{array}\right).
\end{equation}
C'est-à-dire après inversion de la matrice carrée dont le 
 déterminant est $\{\bar{\lambda},\lambda\}$,
\begin{equation}\label{t.22}
\frac{1}{\{\bar{\lambda},\lambda\}}
\left(
\begin{array}{ccc}
\bar{\lambda}'_{\xi}& -\lambda'_{\xi}\\
-\bar{\lambda}'_{x}&\lambda'_{x}
\end{array}\right)
\left(\begin{array}{c}
1\\
0
\end{array}\right)=
\left(\begin{array}{c}
\partial_{z}x_{+}\\
\partial_{z}\xi_{+}
\end{array}\right).
\end{equation}
Nous en tirons les relations suivantes,
\begin{equation}\label{t.23}
\partial_{z}x_{+}=\frac{1}{\{\bar{\lambda},\lambda\}}\bar{\lambda}'_{\xi},\quad
\partial_{\bar{z}}x_{+}=-\frac{1}{\{\bar{\lambda},\lambda\}}\lambda'_{\xi},
\end{equation}
et 
\begin{equation}\label{t.24}
\partial_{z}\xi_{+}=-\frac{1}{\{\bar{\lambda},\lambda\}}\bar{\lambda}'_{x},\quad
\partial_{\bar{z}}\xi_{+}=\frac{1}{\{\bar{\lambda},\lambda\}}\lambda'_{x}.
\end{equation}
De l'équation $\lambda(x,\varphi'_{x}(x,z))-z=0,$ il vient
\begin{equation}
\lambda'_{\xi}(x,\varphi'_{x}(x,z))\varphi''_{x,\bar{z}}(x,z)=0,
\mbox{ soit }\varphi''_{x,\bar{z}}=0,
\end{equation}
et 
\begin{equation}
0=\lambda'_{x}(x,\varphi'_{x})+\lambda'_{\xi}(x,\varphi_{x}')\varphi''_{x,x},
\mbox{ soit }
\varphi''_{x,x}=-\frac{\lambda'_{x}}{\lambda'_{\xi}}(x,\varphi'_{x}).
\end{equation}

Nous sommes maintenant en mesure de démontrer l'égalité (\ref{t.21}).
Nous avons facilement
\begin{align*}
M(z):=\mathrm{Im }\;\frac{\partial}{\partial z}\bigg(\overline{(\partial_{z}
\varphi )(x_{+}(z),z)}\bigg)(z)&=
-\mathrm{Im }\;\frac{\partial}{ \partial{\bar{z}}}\bigg((\partial_{z}
\varphi)(x_{+}(z),z)\bigg)(z)\\
&=-\mathrm{Im}\,(
\varphi''_{x,z}(x_{+},z)\,\partial_{\bar{z}}x_{+} +
\varphi''_{z,\bar{z}}(x_{+},z))
\end{align*}
et, aussi sous l'hypothèse $\varphi(x_{+}(z),z)=0,$
\begin{align*}
\partial_{\bar{z}}(\varphi(x_{+},z))=&
\varphi'_{x}(x_{+},z) \,\partial_{\bar{z}} x_{+}+\varphi'_{\bar{z}}(x_{+},z)=0,\\
\partial_{z}\partial_{\bar{z}}(\varphi(x_{+},z))
=&\varphi'_{x}(x_{+},z)\,\partial_{z}\partial_{\bar{z}} x_{+}+
\varphi''_{x,z}(x_{+},z)\,\partial_{\bar{z}}x_{+}+
\varphi''_{x,x}(x_{+},z)\,|\partial_{z}x_{+}|^{2}\\\
&+\varphi''_{x,\bar{z}}(x_{+},z)\,\partial_{z}x_{+}+\varphi''_{z,\bar{z}}(x_{+},z)=0.
\end{align*}

Le terme $\varphi'_{x}(x_{+},z)\partial_{z}\partial_{\bar{z}} x_{+}$ est réel.
En regroupant tout ce qui précède, nous obtenons donc
\begin{align*}
M(z)= &\mbox{Im }(\varphi''_{x,\bar{z}}(x_{+},z)\,\partial_{z}x_{+})+
|\partial_{z}x_{+}|^{2}\mbox{ Im }\varphi''_{x,x}(x_{+},z)\\
=&
-\frac{|\lambda'_{\xi}|^{2}}{|\{\lambda,\bar{\lambda}\}|^{2}}\mbox{Im }\left(\frac{\lambda'_{x}}{\lambda'_{\xi}}\right)
=
-\frac{1}{|\{\lambda,\bar{\lambda}\}|^{2}} \mbox{Im }\lambda'_{x}\bar{\lambda}'_{\xi}\\
=&
-\frac{1}{|\{\lambda,\bar{\lambda}\}|^{2}} \frac{\{\bar{\lambda},\lambda\} }{2i}
=
-\frac{\{\lambda,\bar{\lambda}\}}{\{\lambda,\bar{\lambda}\}^{2}\,2i} 
=\frac{-1}{2i\{\lambda,\bar{\lambda}\}},
\end{align*}
car $|\{\lambda,\bar{\lambda}\}|^{2}=-\{\lambda,\bar{\lambda}\}^{2}.$
Le côté gauche de (\ref{t.21}) $=\frac{-2}{i\{\lambda,\bar{\lambda}\}}d\mathrm{Re}z
\wedge d\mathrm{Im}z.$
Pour finir, un calcul direct, au moyen de (\ref{t.23}) et (\ref{t.24}), donne
$d\xi_{+}\wedge dx_{+}=\frac{1}{\{\lambda,\bar{\lambda}\} }dz \wedge d\bar{z}.$
La preuve se termine avec l'égalité suivante
\begin{equation}
-\frac{1}{2i}dz\wedge d\bar{z}=d\mathrm{Re}\,z\wedge d\mathrm{Im}\,z.
\end{equation}
\hfill $\square\square$ \medskip

Pour le problème perturbé, nous avons 
\begin{align}\label{pp.35}
\partial_{\bar{z}}\det E_{-+}^{\delta}
&=-\mbox{tr}((\partial_{\bar{z}}R_{+})E_{+}^{\delta}+E_{-}^{\delta}(\partial_{\bar{z}}R_{-}))
\det E_{-+}^{\delta}\\
&=:-k^{\delta}(z)\det E_{-+}^{\delta}.\nonumber
\end{align}
Alors, grâce à la proposition \ref{pp5}, et au fait que $\|\partial_{z}e_{+}\|,\|\partial_{\bar{z}}e_{-}\|
=\mathcal{O}(1/h),$ nous obtenons 
\[ |k^{\delta}-k^{0}|=\mathcal{O}(\frac{\delta \ln (h^{-1})}{h^{3/2}}).\]

On suppose maintenant que $\delta<\frac{h^{3/2}}{\ln (h^{-1})}.$
\begin{prop}
Soit $l^{\delta}$ la  solution de l'équation $\frac{1}{h}\partial_{z}l^{\delta}=k^{\delta},$
donnée par
\[l^{\delta}(z)=\frac{h}{\pi}\int_{\widetilde{W}(z_{0})}\frac{k^{\delta}
(z)}{z-z'}
d\mathrm{Re}z'd\mathrm{Im}z',
\]
Alors $e^{l^{\delta}/h}E_{-+}^{\delta}$ est holomorphe et
$|l^{\delta}-l^{0}|=\ln (h^{-1})\,\mathcal{O}(\frac{\delta}{\sqrt{h}}).$
\end{prop}
\textbf{Preuve.} Nous avons :
\begin{align} 
|(l^{\delta}-l^{0})(z)| &=\left| \frac{h}{\pi}
\int_{\widetilde{W}(z_{0})}\frac{(k^{0}-k^{\delta})(z')}{z-z'}
d\mbox{Re} z'd\mbox{Im}z' \right|\nonumber\\
&\le h\|k^{\delta}-k^{0}\|_{L^{\infty}(\widetilde{W}(z_{0}))}\int_{\widetilde{W}(z_{0})}
 \frac{1}{|z-z'|}d\mbox{Re} z'd\mbox{Im}z'\nonumber\\
&= \mathcal{O}(\frac{\delta\ln (h^{-1})}{\sqrt{h}}).\label{pp.37}
\end{align}
\hfill $\square$ \medskip
\section{Estimation de la probabilité que $\det E_{-+}^{\delta }$ soit petit}
Soit $z_{0}$ un point appartenant à $\Omega\Subset\Lambda(p).$ 
L'ensemble $\Omega$ a été introduit dans l'hypothèse \ref{i3}. Rappelons également
que $\widetilde{W}(z_{0})$ est l'ensemble intervenant dans le théorème 
\ref{gg6} et  que $\delta$ vérifie
\[\delta\ll\frac{h^{3/2}}{|\ln h|}.\]

En nous restreignant   à $\|Q\|\le  |\ln h|,$ nous savons  alors que 
$E^{\delta}_{-+}$ s'écrit
\begin{equation}\label{w.6}
E_{-+}^{\delta}=E_{-+}+ E_{-}Q^{\delta}E_{+}+
E_{-}Q^{\delta}\left(\sum_{k\ge 0}E(Q^{\delta}E)^{k}\right)Q^{\delta}E_{+},
\end{equation}
où $Q^{\delta}:=\delta Q.$
Le terme entre parenthèses, désigné par $\widetilde{E}$, 
est $\mathcal{O}(1/\sqrt{h}).$ 
Nous trouvons alors
\begin{align}
\det E_{-+}^{\delta}=&\sum_{\pi\in S_{\beta}}
\prod_{1\le\nu\le\beta }
(\mbox{sign}(\pi))
\left(
\langle Q^{\delta}e^{\pi (\nu)}_{+},e^{\nu}_{-}\rangle+
\langle Q^{\delta}\widetilde{E}Q^{\delta}e_{+}^{\pi(\nu)},e^{\nu}_{-}\rangle
\right)\nonumber\\
&+\mathcal{O}(h^{\infty}).
\end{align}

L'opérateur $\widetilde{E}$  satisfait la condition 
de non-propagation 
du support, au sens
défini dans le théorème \ref{gg6}. Gr\^ace à l'hypothèse \ref{i3}
(et plus particulièrement, au fait que $x_{+}^{\nu}\ne x_{-}^{\kappa}$ 
pour $\nu\ne\kappa$), alors
\begin{equation}\label{w.8}
\det E_{-+}^{\delta}=
\prod_{1\le \nu\le \beta }\left(
\langle Q^{\delta} e^{\nu}_{+},e^{\nu}_{-}\rangle+
\langle Q^{\delta}\widetilde{E}Q^{\delta}e_{+}^{\nu},e^{\nu}_{-}\rangle
\right)+\mathcal{O}(h^{\infty}).
\end{equation}
 Remarquons ensuite, que si
\begin{equation}\label{w.9}
|\langle Q^{\delta}e^{\nu}_{+},e^{\nu}_{-}\rangle | \ge \frac{3}{2} x^{\frac{1}{\beta}},
\mbox{ et }
|\langle Q^{\delta}\widetilde{E}Q^{\delta}e_{+}^{\nu},e^{\nu}_{-}\rangle|
\le \frac{1}{2} x^{\frac{1}{\beta}},
\end{equation}
nous avons
\begin{equation}
|\langle Q^{\delta} e^{\nu}_{+},e^{\nu}_{-}\rangle
+
\langle Q^{\delta}\widetilde{E}Q^{\delta}e_{+}^{\nu},e^{\nu}_{-}\rangle|
\ge x^{\frac{1}{\beta}}.
\end{equation}
Ce qui entra\^ine la minoration suivante (confère l'inégalité 
$\mathbb{P}(A\cap B)\ge  \mathbb{P}(A)+\mathbb{P}(B)-1,$
pour deux évènements $A$ et $B$ quelconques) si $x>0$,
%\textit{$x$ est supposé minoré par une certaine puissance de $h$} 
\begin{multline}\label{w.11}
\mathbb{P}(|\det E_{-+}^{\delta}|\ge x+\mathcal{O}(h^{\infty}))\ge
\mathbb{P}(\|Q\|\le |\ln h|)-2\beta\\
+\sum_{1\le\nu\le \beta} 
(\mathbb{P}(
|\langle Q^{\delta}e_{+}^{\nu},e_{-}^{\nu}\rangle|
\ge \frac{3}{2} x^{\frac{1}{\beta}})+
\mathbb{P}(
|\langle Q^{\delta}\widetilde{E}Q^{\delta}e_{+}^{\nu},e^{\nu}_{-}\rangle|
\le\frac{1}{2} x^{\frac{1}{\beta}})).
\end{multline}

Nous avons une estimation de  $\mathbb{P}(\|Q\|\le |\ln h|)$
(corollaire \ref{pp4}).
Il nous reste donc à estimer les probabilités des deux derniers 
termes de (\ref{w.11}). %ce que nous ferons par la suite. 
Commen\c cons 
par étudier 
$
\mathbb{P}(
|\langle Q^{\delta}e_{+}^{\nu},e_{-}^{\nu}\rangle|
\ge \frac{3}{2} x^{\frac{1}{\beta}}),
$
puis
$
\mathbb{P}(
|\langle Q^{\delta}\widetilde{E}Q^{\delta}e_{+}^{\nu},e^{\nu}_{-}\rangle|
\le\frac{1}{2} x^{\frac{1}{\beta}}).
$
\begin{lemme}
Il existe $C>0$ tel que pour tout $\nu\in1,\ldots,\beta$
\begin{equation}\label{w.13}
\mathbb{P}( \|Q^{\delta}e_{+}^{\nu}\|\le x)\ge 1-C\exp(
-\frac{1}{C\delta^{2}}\,x^{2}).
\end{equation}
La conclusion est la même pour $\|\delta Q^{\ast}e_{-}^{\nu}\|.$ 
Soit
\begin{equation}\label{w.14}
\mathbb{P}(
\langle Q^{\delta}\widetilde{E}Q^{\delta}e_{+}^{\nu},e^{\nu}_{-}\rangle
\le x)
\ge 1-\widetilde{C}\exp(
-\frac{1}{\widetilde{C}\delta^{2}}\,x\sqrt{h}). \nonumber
\end{equation}
\end{lemme}
\textbf{Preuve.}
Avec l'aide de (\ref{pp.4}), nous trouvons que
\begin{equation}
\mathbb{P}( \|Q\| \le \frac{x}{\|e_{+}^{\nu}\|})
\ge 1-C\exp(-\frac{1}{C}x^{2}).
\end{equation}
pour une constante $C>0.$ 
Puisque 
$\|\widetilde{E}\| \lesssim h^{-\frac{1}{2}}$, nous avons
\[
|\langle Q\widetilde{E}Q e_{+}^{\nu},e^{\nu}_{-}\rangle|=
|\langle \widetilde{E}\,Qe_{+}^{\nu},Q^{\ast}e_{-}^{\nu}\rangle|
\lesssim 
\|Qe_{+}^{\nu}\|\,\|Q^{\ast}e_{-}^{\nu}\|\,h^{-\frac{1}{2}}.
\]
Terminons en appliquant  l'inégalité citée avant (\ref{w.11}).
\hfill $\square$ \medskip

Cherchons maintenant à préciser la loi de probabilité de 
$\langle Qe_{+},e_{-}\rangle,$ où $e_{\pm}\equiv e_{\pm}^{\nu}$ 
pour un $\nu$ fixé. Un calcul direct 
de $\langle Qe_{+},e_{-}\rangle$ donne
\begin{equation}\label{ss.10}
\sum_{\alpha,i,j}\langle Q_{\alpha}^{i,j}(x) (hD_{x})^{\alpha}e_{+,j},
e_{-,i}\rangle\\
=\sum_{\alpha,i,j,k} q_{\alpha,k}^{i,j}\,
\langle e_{k}\,
(hD_{x})^{\alpha}e_{+,j}, e_{-,i}\rangle.
\end{equation}
où $e_{\pm,i}$ sont les  coordonnées de $e_{\pm}$ et $e_{k}:=\frac{e^{ikx}}{\sqrt{2\pi}}.$ 
Il devient alors évident que $\langle Qe_{+},e_{-}\rangle$ suit la loi gaussienne
$\mathcal{N}(0,\sigma^{2}),$ la variance satisfaisant 
\begin{equation}\label{ss.11}
\sigma^{2}(h)=\sum_{\alpha,i,j,k}(\sigma_{\alpha,k}^{i,j}(h))^{2}\;
|\langle e_{k}\,(hD_{x})^{\alpha} e_{+,j},e_{-,i}\rangle|^{2}.
\end{equation}

Reste alors à donner le comportement de $\sigma^{2}(h)$
lorsque $h\to 0.$  
\begin{lemme}\label{ss3}
Il existe $C>0$ tel que pour tout $z$ dans $\widetilde{W}(z_{0}),$\\
(a) si $|k|\le\frac{1}{hC}$  on a 
$|\langle e_{k}\,(hD_{x})^{\alpha}e_{+,j},e_{-,i}\rangle |=\mathcal{O}(h^{\infty}),$\\
(b) si $\frac{1}{hC}\le|k|\le\frac{C}{h}$ on a 
$|\langle e_{k}\,(hD_{x})^{\alpha}e_{+,j},e_{-,i}\rangle |=\mathcal{O}(1),$\\
(c) si $|k|\ge\frac{C}{h}$  on a 
$|\langle e_{k}\,(hD_{x})^{\alpha}e_{+,j},e_{-,i}\rangle |=\mathcal{O}(1/|k|^{\infty}).$
\end{lemme}
\textbf{Preuve.} 
Pour $(b),$ utilisons l'inégalité de Cauchy-Schwarz
\[|\langle e_{k}\,(hD_{x})^{\alpha}e_{+,j},e_{-,i}\rangle |\le
\|e_{-,i}\|\,\|(hD_{x})^{\alpha} e_{+,j}\|=\mathcal{O}(1).\]

Pour $(a)$ et $(c),$ remarquons d'abord que 
$\langle e_{k}\,(hD_{x})^{\alpha}e_{+,j},e_{-,i}\rangle $ est une intégrale du type
\[
h^{-1/2}\int_{0}^{2\pi}e^{-\frac{i}{h}\phi(x,z)+ikx}
a(x;h)dx,
\]
où $\phi=\overline{\varphi_{-}}-\varphi_{+}$ satisfait  
$\phi'(x(z))=\xi_{-}(z)-\xi_{+}(z)\ne 0$ si on note
$x(z):=x_{+}(z)=x_{-}(z)$, et  $a$
est un symbole de classe $S(1)$ à support compact, 
qui contient $x(z).$ Ecrivons
$\varphi:= -\frac{\phi}{h}+kx.$ Il existe $C>0$ pour lequel 
\[\forall z\in \widetilde{W}(z_{0}),\quad \forall |k|\notin [\frac{1}{hC},
\frac{C}{h}],\quad
|\varphi'_{x}(x,z)|\ge \frac{1}{C}\max (|k|,1/h).\]
Nous nous servons du fait que $\inf |\phi'|\ne 0$ dans un voisinage de $x(z)$
avant de procéder par intégration par partie pour trouver
\begin{align}
\langle e_{k}&\,(hD_{x})^{\alpha}e_{+,j},e_{-,i}\rangle =
\frac{1}{i} \int e^{i\varphi(x)}a_{n}(x)dx,
\nonumber
\\
&a_{n}:=\left(-\frac{d}{dx}\circ \frac{1}{\varphi'}\right)^{n}(a)
=\mathcal{O}\left(
(\min(1/|k|,h))^{n}\right).
\end{align}
\hfill $\square$ \medskip
\begin{prop}\label{ss3bis} Soit $Q$ vérifiant l'hypothèse \ref{i4}.
Il existe $\widetilde{C}>0$  tel que nous avons 
$\langle Qe_{+},e_{-}\rangle\sim \mathcal{N}(0,\sigma^{2}),$
où la variance vérifie 
\[\frac{1}{\widetilde{C}} h ^{2\rho-1/2}\le \sigma^{2}(h)\le 
\widetilde{C}h^{2\rho-1/2}.\]
\end{prop}
\textbf{Preuve.}
Pour la borne inférieure, il suffit de montrer  qu'il existe $i$ et $j$
pour lesquels nous avons pour  $\alpha=\alpha_{1},$
\begin{equation}\label{ss.12}
h^{2\rho-1/2}\lesssim\sum_{k\in\Z}(\sigma_{\alpha_{1},k}^{i,j})^{2}\;
|\langle e_{k}\,(hD_{x})^{\alpha_{1}}  e_{+,j},e_{-,i}\rangle|^{2}
\quad(\le \sigma^{2}(h)).
\end{equation}
Considérons alors pour chaque $i,j,$ la somme (\ref{ss.12}). 
Si nous découpons la  sommation sur $k$ en trois, suivant $k\gg 1/h,$ 
$k\ll 1/h$ et $k\sim 1/h,$ 
 (\ref{ss.12})  s'écrit alors de fa\c con précise
\begin{multline}\label{ss.14}
\sum_{|k|<\frac{1}{hC},\,
|k|>\frac{C}{h}}(\sigma_{\alpha_{1},k}^{i,j})^{2}\;
|\langle e_{k}\,(hD_{x})^{\alpha_{1}}  e_{+,j},e_{-,i}\rangle|^{2}\\
+
\sum_{\frac{1}{Ch}\le|k|\le \frac{C}{h}}(\sigma_{\alpha_{1},k}^{i,j})^{2}\;
|\langle e_{k}\,(hD_{x})^{\alpha_{1}}  e_{+,j},e_{-,i}\rangle|^{2},
\end{multline}
où $C>0$ est la constante du lemme \ref{ss3}. Celle-ci montre que le premier
terme est $\mathcal{O}(h^{\infty}).$
Puis, gr\^ace à l'identité de Parseval et l'hypothèse \ref{i4} sur la minoration 
des variances $\sigma^{i,j}_{\alpha_{1},k}$,
le second terme de (\ref{ss.14}) est  pour tout $i$ et $j$ 
\begin{align}\label{ss.15}
\gtrsim & h^{2\rho}
\bigg(\|((hD_{x})^{\alpha_{1}}e_{+,j})\,\overline{e_{-,i}}\|^{2}
-\sum_{|k|<\frac{1}{Ch},\,|k|>\frac{C}{h}}
|\langle e_{k} \,(hD_{x})^{\alpha_{1}} e_{+,j},e_{-,i}\rangle |^{2}
\bigg)
\nonumber\\
\gtrsim & h^{2\rho}
\big(\|((hD_{x})^{\alpha_{1}}e_{+,j})\,\overline{e_{-,i}}\|^{2}
+\mathcal{O}(h^{\infty})\big).
\end{align}

Observons d'une part que, puisque $\|e_{\pm}\|= 1,$  
certaines coordonnées de $e_{\pm}$
sont elliptiques, c'est-à-dire de terme principal non nul; 
et d'autre part que, 
si $e_{+,j}$ et $e_{-,i}$ sont elliptiques, alors nous avons 
\begin{equation}\label{ss.16}
\|((hD_{x})^{\alpha_{1}}e_{+,j})\overline{e_{-,i}}\|^{2}\asymp h^{-1/2}\;
(\mbox{car } \xi_{+}\ne 0).
\end{equation}
Dans le cas où $\xi_{+}=0,$ la relation ci-dessous ne tient plus. En effet,
le symbole principal de $((hD_{x})^{\alpha_{1}}e_{+,j})$ s'annule au point critique.

Par conséquent, si $e_{+,j}$ et $e_{-,i}$ sont choisis elliptiques dans (\ref{ss.12}), 
nous avons le résultat demandé.  
Il suffit d'appliquer la  même procédure pour montrer la borne supérieure.
\hfill $\square$ \medskip

Si  $X$ suit la loi gaussienne complexe $\mathcal{N}(0,\sigma^{2})$ alors 
\begin{equation}\label{ss.21}
\mathbb{P}(|X|\ge x)=\exp(-\frac{x^{2}}{\sigma^{2}}).
\end{equation}

En résumé: si $z\in\widetilde{W}(z_{0})$ et $x>0$ alors nous avons l'estimation 
suivante
du membre de gauche de  (\ref{w.11}) 
\begin{align}\label{ss.22}
\mathbb{P}(\det E^{\delta}_{-+}(z)&>x+\mathcal{O}(h^{\infty}))
\ge 1-Ce^{-\frac{1}{C}(\ln h)^{2}}
\nonumber\\
&+\beta(\exp [-C\frac{x^{2/\beta}}{\delta^{2}h^{2\rho-1/2}} ]-1)
-C\beta \exp[ -\frac{x^{1/\beta}}{C\delta^{2} h^{-1/2}} ].
\end{align}
pour une constante $C>0.$
En prenant $x$ de l'ordre de $h^{(\rho+\epsilon-1/4)\beta}\delta^{\beta},$ 
avec $\delta$ minoré par une puissance de $h$, pour que $x$ soit le terme
dominant de $x+\mathcal{O}(h^{\infty}),$
nous pouvons proposer:
\begin{prop} Pour tous $z\in\widetilde{W}(z_{0}),$ $\epsilon>0,N_{0}\gg1,$
et 
\[h^{N_{0}}\ll\delta\ll h^{\rho+\epsilon+\frac{1}{4}}|\ln h|^{-2}\]
nous avons 
\begin{align}\label{ss.23}
\mathbb{P}(|\det E^{\delta}_{-+}(z)|\ge h^{(\rho+\epsilon-1/4)\beta}
\delta^{\beta})\ge 1-Ch^{2\epsilon},
\end{align}
pour une constante $C>0.$
\end{prop}
\section{Preuve du Théorème \ref{i6}}
Grâce à $l^{\delta}=l^{0}+\mathcal{O}(\frac{\delta|\ln h|}{\sqrt{h}}),$ 
de (\ref{w.8}) nous obtenons  
que
\begin{align}
| e^{\frac{l^{\delta}}{h}}&\det E^{\delta}_{-+}|\le
e^{\frac{\mathrm{Re}\, l^{0}}{h}+\mathcal{O}(\frac{\delta |\ln h|}{h^{3/2}})}\nonumber\\
&\times \prod_{1\le\nu\le\beta}
\left(
\big(\|\delta Q\|
\|e_{+}^{\nu}\|
\|e_{-}^{\nu}\|+
\|\delta Q\|^{2}
\|e_{+}^{\nu}\|
\|e_{-}^{\nu}\|
\mathcal{O}(\frac{1}{\sqrt{h}})\big)
+\mathcal{O}(h^{\infty})\right).
\end{align}
Puisque $\|Q\|\ll |\ln h|$ et $\delta\ll h^{\rho+\epsilon+\frac{1}{4}}|\ln h|^{-2},$
il en résulte que
\begin{equation}\label{ss.25}
| e^{\frac{l^{\delta}}{h}}\det E^{\delta}_{-+}(z)|\le e^{\frac{1}{h}\mathrm{Re}\,l^{0}(z)}, 
\quad\forall z\in \widetilde{W}(z_{0}).
\end{equation}
Introduisons la fonction holomorphe 
\begin{equation}\label{ss.26}
 F_{\delta}(z,h):=e^{\frac{l^{\delta}}{h}}\det E_{-+}^{\delta}(z),
 \quad z\in \widetilde{W}(z_{0}).
\end{equation}
\begin{corol}\label{ss4}
 Soient $z_{0}$ un point de $\Omega\Subset\Lambda(p)$ et
 $\epsilon,N_{0}>0.$ Il existe un voisinage de $z_{0}$ noté $\widetilde{W}(z_{0})$ 
inclus dans  $\Omega$ tel que, si 
$h^{N_{0}}\ll\delta\ll h^{\rho+\epsilon+\frac{1}{4}}|\ln h|^{-2},$
alors il existe $C,\widetilde{C}>0$ telles que\\
(a) avec une probabilité $\ge 1-Ce^{-\frac{1}{C}(\ln h)^{2}}$
nous avons
\begin{equation}\label{ss.26bis}
\ln | F_{\delta}(z,h)|\le \frac{1}{h}\mathrm{Re}\,l^{0}(z),
\end{equation}
pour tous les $z$ dans $\widetilde{W}(z_{0}).$\\
(b) pour chaque $z$ de $\widetilde{W}(z_{0}),$ $\epsilon>1/4$ nous avons 
\begin{align}
 \ln | F_{\delta}(z,h)|&\ge 
 \frac{1}{h}(\mathrm{Re}\,l^{0}(z)-C\frac{\delta |\ln h|}{\sqrt{h}}-
 h \beta|\ln (h^{\rho+\epsilon-1/4}\delta)|)\nonumber\\
 &\ge 
 \frac{1}{h}(\mathrm{Re}\,l^{0}(z)-h(\widetilde{C}+ \beta)
 |\ln (h^{\rho+\epsilon-1/4}\delta)|),
 \label{ss.27}
 \end{align}
avec une probabilité $\ge 1-Ch^{2\epsilon}$.
\end{corol}
\textbf{Preuve.}
$(a)$ découle de (\ref{ss.25}) et du corollaire \ref{pp4},
pour $(b)$ il faut se référer à (\ref{ss.23}).
\hfill $\square$ \medskip

Nous pouvons maintenant répéter les arguments de \cite{Ha1, Ha2}. 
Rappelons une proposition  de \cite{Ha1}, qui reste valable 
pour des contours $C^{2}$ par morceaux (confère
théorème \ref{i6}).
En effet, la même preuve permet d'avoir une frontière $C^{2}$ avec un nombre fini
de points anguleux.
\begin{prop}\label{ss5}
Soient $\Omega\Subset\C,$ $\Gamma\subset\Omega$ un domaine
à bord $C^{2}$ par morceaux et $\phi\in C^{\infty}(\Omega,\R).$ 
Soit $f$ une fonction holomorphe dans $\Omega$ vérifiant
\begin{equation}\label{ss.28}
|f(z,h)|\le e^{\phi(z)/h},\quad z\in\Omega.
\end{equation}
Supposons qu'il existe $\tilde{\epsilon}\ll 1, z_{k}\in \Omega,
k\in K$ tels que 
\begin{align}
&\partial\Gamma\subset\bigcup_{k\in K}D(z_{k},\sqrt{\tilde{\epsilon}}),
\quad \#K=\mathcal{O}
(\frac{1}{\sqrt{\tilde{\epsilon}}}),\\
 & |f(z_{k},h)|\ge e^{\frac{1}{h}(\phi(z_{k})-\tilde{\epsilon})},
 \quad k\in K,\label{ss.30}
 \end{align}
alors 
\[
\#(f^{-1}(0)\cap \Gamma)=
\frac{1}{2\pi h}
\iint_{\Gamma}\Delta \phi\,
d(\mathrm{Re} z) d(\mathrm{Im} z)
+ \mathcal{O}(\frac{\sqrt{\tilde{\epsilon}}}{h}).
\]
\end{prop}

Nous pouvons  appliquer la proposition 
avec $\tilde{\epsilon}=h(\widetilde{C}+\beta)|\ln (h^{\rho+\epsilon-1/4}\delta)|,$
$\phi=\mbox{Re}\,l^{0}$
et $f=F_{\delta}.$ L'évènement 
(\ref{ss.28}) a la même probabilité de se réaliser que 
l'évènement (\ref{ss.26bis}). L'évènement (\ref{ss.30}) 
se réalise avec 
une probabilité
\begin{align}
&\ge1-Ch^{2\epsilon}\,(\#K)\nonumber\\
&\ge 1-\tilde{C}\frac{h^{2\epsilon}}
{\sqrt{h|\ln (h^{\rho+\epsilon-1/4}\,\delta)|}}.
\end{align} 

Compte tenu de la remarque après la  proposition \ref{pp6}, 
nous sommes maintenant en possession du résultat suivant $(\gamma_{1}+
\frac{1}{4}=\epsilon,$ $\gamma_{1}>0)$:
\begin{theo} \label{ss6}
Soient $z_{0}$ un point de $\Omega\Subset\Lambda(p),$ et $N_{0}\gg 1.$ 
 Il existe un voisinage 
$\widetilde{W}(z_{0})$ de $z_{0},$ tel que si $\Gamma$ est un ouvert relativement compact
dans $\widetilde{W}(z_{0}),$   à bord $C^{2}$ par morceaux,
et $\gamma_{1}>0,$  alors il existe $C>0$ tel que si
\begin{align}\label{ss.31}
h^{N_{0}}\ll\delta\ll h^{\rho+\gamma_{1}+\frac{1}{2}}
|\ln h|^{-2},
\end{align}
alors le spectre de $P-\delta Q_{\omega}$ est discret, et,
le nombre $N(P-\delta Q_{\omega},\Gamma)$  de valeurs propres de 
$P-\delta Q_{\omega}$
dans $\Gamma$ satisfait
\begin{align}\label{ss.31bis}
|N(P-\delta Q_{\omega},\Gamma)-\frac{1}{2\pi h} \iint m_{\Gamma}\, dxd\xi |\le C
h^{-\frac{1}{2}}|\ln(h\delta)|^{\frac{1}{2}}
\end{align}
avec une probabilité 
\[\ge 1-C
h^{2\gamma_{1}}|\ln(h\delta)|^{-\frac{1}{2}}.\]
Rappelons que $m_{\Gamma}$ a été introduit en (\ref{i.15}).
\end{theo}

Soit $\Gamma\Subset\Omega.$
En  recouvrant $\Gamma$ par un nombre fini de 
$\Gamma_{k}\subset \widetilde{W}(z_{k}),$ 
nous obtenons le théorème \ref{i6}.

Précisons que nous n'avons aucune hypothèse garantissant 
que $\Sigma(p)\ne\C.$ Le spectre discret
est une conséquence de la perturbation aléatoire.
Ce fait résulte de la théorie de Fredholm analytique
(impliquant que le spectre d'un opérateur de Fredholm, 
d'indice zéro, est soit discret soit $\C$) et de l'estimation probabiliste 
(\ref{ss.23}) ($\det E^{\delta}_{-+}(z)$ ne s'annule pas 
avec une forte probabilité).

Nous allons maintenant donner un résultat similaire 
concernant l'asymptotique de Weyl
pour une famille de domaines. Rappelons d'abord la proposition suivante 
qui est le cas uniforme de la proposition \ref{ss5}.
\begin{prop}\label{ss8}
Soit un domaine $\Omega\Subset\C$, et 
$\mathcal{G}$ une famille de domaines inclus dans $\Omega$. 
Soit $C_{0}>0$ une constante indépendan\-te de 
$\mathcal{G}$.
Supposons que 
\begin{equation}\label{ss.32}
\forall\,\Gamma\in\mathcal{G},\quad
\partial \Gamma=\bigcup_{j=1}^{N}\gamma_{j}([a_{j},b_{j}]),\quad N\le C_{0},
\end{equation}
où $\gamma_{j}:[a_{j},b_{j}]\to\C$ est $C^{2}$ 
%($\gamma_{j}$ dépend évidemment de $\Gamma$) 
avec
\begin{align}
&\forall j,\quad 0<a_{j}<b_{j}\le C_{0},\\
&\forall j,\quad \frac{1}{C_{0}}\le |\dot{\gamma}_{j}(t)|\le C_{0},\quad
|\ddot{\gamma}_{j}(t)|\le C_{0},\\
&\gamma_{j}(b_{j})=\gamma_{j+1}(a_{j+1}),\;j\in \Z/N\Z. \label{ss.33}
\end{align}
Soit $\phi\in C^{\infty}(\Omega,\R)$ et  
$f$ une fonction holomorphe dans $\Omega$ avec
\begin{equation}\label{ss.34}
|f(z;h)|\le e^{\frac{\phi(z)}{h}},\;\forall z\in\Omega.
\end{equation}
Pour un réseau carré de points $z_{k}\in\Omega$ de maille 
$\frac{\sqrt{\tilde{\epsilon}}}{2}$, 
$0<\tilde{\epsilon}\ll 1$ avec
\[\Omega\subset\bigcup_{k\in K}D(z_{k},\frac{\sqrt{\tilde{\epsilon}}}{2}),
\quad |K|\le \frac{C}{\tilde{\epsilon}}\]
et \begin{equation}\label{ss.35}
|f(z_{k};h)|> e^{\frac{1}{h}(\phi(z_{k})-\tilde{\epsilon})},
\end{equation}
alors 
\begin{align*}
&\exists \,D>0,\;\forall \,\Gamma\in\mathcal{G},\\
&|\#(f^{-1}(0)\cap \Gamma )-\frac{1}{2\pi h}\iint_{\Gamma}\Delta\phi\;
d(\mathrm{Re}\;z)d(\mathrm{Im}\;z)|\le
D\frac{\sqrt{\tilde{\epsilon}}}{h}.
\end{align*}
\end{prop}

A la section 6.3 (``Preuve du Théorème 1.9'') de \cite{Ha2}, M.~Hager  démontre la proposition \ref{ss8} 
pour une famille différente de $\Gamma.$ Pour la preuve, nous suivrons ici
la démarche de Hager de la section 6.3 (du début de la démonstration à (6.16))
conjuguée avec le lemme suivant :
\begin{lemme}\label{an2}
Soit une famille de lacets simples $\gamma_{j}$ ($j\in J$) dans $\C$ de
classe $C^{2}$.  Paramétrisons les lacets
$\gamma_{j}:\,[0,1]\owns t \to \gamma_{j}(t)\in\C.$ 
Supposons également qu'il existe $C_{0}>0$ tel que
\[\forall j\in J,\;\forall t\in[0,1],\quad 
\frac{1}{C_{0}}\le |\dot{\gamma}_{j}(t)|\le C_{0},
\quad |\ddot{\gamma}_{j}(t)|\le C_{0}.\]
%Soit $z\in \C$ un point satisfaisant $d(z,\gamma_{i})\le r.$
Alors il existe une constante $C>0,$ indépendante de $j\in J,$ telle que 
pour tout $r\ll1/C_{0}^{3}$ chaque composante connexe de $\gamma_{j}\cap D(z,r)$ est de longueur
$\le Cr,$ où le point $z$ est donné. 
\end{lemme}

%Le lemme nous apprend que, dans un disque donné, la longueur d'un morceau de courbe 
 %de la famille $J$ est limitée.
%Pour la démonstration \ref{ss8} nous avons besoin en effet d'une 
%majoration uniforme (voir \cite{Ha2}).

\noindent
\textbf{Preuve.}
Posons $f_{i}(t):=\frac{1}{2}|\gamma_{i}(t)-z|^{2}.$
Un calcul montre (en omettant les indices) que
\begin{eqnarray*}
\dot{f}(t)&=&\langle\gamma(t)-z,\,\dot{\gamma}(t)\rangle,\\
\ddot{f}(t)&=&|\dot{\gamma}(t)|^{2} + \langle\gamma(t)-z,\ddot{\gamma}(t)\rangle,
\end{eqnarray*}
où $\langle ., .\rangle$ représente le produit scalaire dans $\R^{2}$.

Pour la suite nous travaillons dans $D(z,r)$, c'est à dire que nos $t$
vérifient $\gamma(t)\in D(z,r)$. 
   
Grâce à l'inégalité de Cauchy-Schwarz et les hypothèses du lemme, 
il existe $C_{1}>0$ tel que pour $r\ll 1/C_{0}^{3}$ 
\begin{eqnarray}
\ddot{f}(t)&\ge&|\dot{\gamma}(t)|^{2}- \;|\gamma(t)-z|\;|\ddot{\gamma}(t)| \label{an0}\\
&\ge& \frac{1}{C_{0}^{2}}-rC_{0}\ge \frac{1}{C_{1}}>0.\nonumber
\end{eqnarray}

Dans chaque composante connexe de $\gamma\cap D(z,r)$, il existe un temps $t_{1}$ pour lequel
$f(t)$ est minimum, dès lors $\dot{f}(t_{1})=0$. 
La formule de Taylor avec reste intégrale 
donne alors 
$$f(t)=f(t_{1})+\int_{t_{1}}^{t}(x-t_{1})\ddot{f}(x)dx$$
de là
$$f(t)\ge f(t_{1})+\frac{1}{2C_{1}}(t-t_{1})^{2}.$$
soit $$|t-t_{1}|\le r\sqrt{2C_{1}}.$$

Ainsi la longueur de chaque composante connexe de $\gamma\cap D(z,r)$
est majorée par 
$$ \int _{|t-t_{1}|\le r\sqrt{2C_{1}}}|\dot{\gamma}(t)|dt
\le r\sqrt{8C_{1}C_{0}^{2}}.$$
%(Toutes les estimations sont bien sûr uniformes par rapport à $j$) 
\hfill $\square$ \medskip

La proposition \ref{ss8} nous conduit donc au résultat qui suit :
\begin{theo}\label{ss9}
Soit $\mathcal{G}$ une famille de domaines $\Gamma\Subset\Omega,$
vérifiant les hypothèses de la proposition \ref{ss8}. Nous supposons que l'hypothèse 
\ref{i3} est satisfaite. Soient $\gamma_{2}>0,$ 
$\delta\ll h^{\rho+\gamma_{2}+\frac{3}{4}}|\ln h|^{-2}$  et minoré par une puissance de $h,$
alors avec une probabilité
\[\ge 1-Ch^{2\gamma_{2}}|\ln(h\delta)|,\]
nous avons (\ref{ss.31bis}) avec une constante $C$ indépendante de $\Gamma.$
\end{theo}

\section{Réduction semiclassique}
Intéressons-nous maintenant à la distribution des grandes valeurs propres de $P-Q_{\omega}.$
Rappelons que $P$ et $Q_{\omega}$ s'écrivent respectivement
\[P=\sum_{0\le\alpha\le m}A_{\alpha}(x)D_{x}^{\alpha},
\quad
Q_{\omega}=\sum_{\alpha_{0}\le \alpha\le\alpha_{1}}Q_{\alpha}(x)D_{x}^{\alpha},
\]
 où les entrées de $Q_{\alpha}$ sont des séries de Fourier aléatoires. 
 $Q_{\alpha}$ satisfait l'hypothèse \ref{i4} et  $P$ est elliptique au sens
 classique. Nous allons pour commencer
restreindre le paramètre spectral $z$ au domaine $\Omega_{R},$
où $\Omega_{R}=R\Omega_{1},$ $R\gg 1,$
avec $\Omega_{1}\Subset\Omega\subset\Lambda(p_{m}).$
Puisque $\Lambda(p_{m}),\Omega$ sont des cônes, nous avons 
pour tout $R\ge1,$ $\Omega_{R}\Subset\Omega\subset\Lambda(p_{m}).$

Pour $z\in\Omega_{R},$ nous ramenons l'étude de $P-z$ à un problème 
semiclassique en divisant par $R.$ Nous sommes donc invités à étudier, en posant
$h^{m}R=1,$ l'opérateur
\begin{equation}\label{w.1}
P^{0}-w-\widetilde{Q}_{\omega}=
h^{m}(P-z-Q_{\omega}),\quad w:=\frac{z}{R}\in \Omega_{1}.
\end{equation}
Le symbole principal semiclassique de $P^{0}$ est alors $p_{m}.$
Nous avons 
\begin{equation}
\widetilde{Q}_{\omega}=\sum_{\alpha_{0}\le\alpha\le \alpha_{1}}
h^{m-\alpha}Q_{\alpha}(x)(hD_{x})^{\alpha},
\quad \widetilde{Q}_{\omega}^{0}:=
h^{-(m-\alpha_{1})}\widetilde{Q}_{\omega}.
\end{equation}

Reprenons les notations $\delta$ et $\Lambda(p_{m})$ 
introduites dans le cadre semiclassique.
$P^{0}$  satisfait l'hypothèse d'ellipticité \ref{i1}, et par l'hypothèse \ref{i10},
les points $\rho_{\pm}\in p_{m}^{-1}(w),$ $w\in \Omega_{1}$ vérifient
la condition \ref{i3}. De plus,
la perturbation $\widetilde{Q}_{\omega}^{0}$ 
entre bien dans le cadre de \ref{i4}. 
Dans ces conditions,
Nous pouvons appliquer le théorème \ref{ss6}, à
$P^{0}-\delta\widetilde{Q}^{0}_{\omega},$ avec $\delta=h^{m-\alpha_{1}}.$
La condition (\ref{ss.31}),
$\delta=h^{m-\alpha_{1}}\ll h^{\rho+\gamma_{1}+\frac{1}{2}}|\ln h|^{-2},$
 équivaut ici à
\begin{equation}\label{w.2}
m-\alpha_{1}>\rho+\gamma_{1}+\frac{1}{2}
\end{equation}
et dans le cas d'une famille de domaines $m-\alpha_{1}>\rho+\gamma_{2}+\frac{3}{4}.$

Notons pour tout $(x,\xi)\in T^{\ast}S^{1}$ et $\Gamma\subset\C,$
\begin{equation}\label{w.3}
m_{\Gamma}:=\#(\sigma(p_{m}(x,\xi))\cap \Gamma).
\end{equation}

Nous avons les égalités suivantes :
\begin{align}\label{w.30}
N(P^{0}-\delta\widetilde{Q}^{0}_{\omega},\Gamma)&=N(P-Q_{\omega},R\Gamma)\\
\frac{1}{2\pi h}\iint m_{\Gamma}\,dxd\xi&=
\frac{1}{2\pi }\iint m_{R\Gamma}\,dxd\xi,
\end{align}
Ce qui implique qu'avec une probabilité  
$\ge 1-C R^{-2\gamma_{1}/m}(\sqrt{\ln R})^{-1}$ 
nous avons
\begin{equation}\label{w.32}
|N(P-Q_{\omega},R\Gamma)-\frac{1}{2\pi } \iint m_{R\Gamma}\, dxd\xi | \le C
R^{1/(2m)}\sqrt{\ln R}
\end{equation}

Si $\mathcal{G}$ une famille de domaines $\Gamma\Subset\Omega_{1},$
vérifiant les hypothèses du théorème \ref{ss8}, alors 
avec une probabilité $\ge 1-C R^{-2\gamma_{2}/m}(\sqrt{\ln R})^{-1}$ 
nous avons (\ref{w.32}) avec une constante $C>0$ indépendante de $\Gamma.$

\section{Preuve du Théorème \ref{j8}}
Nous nous intéressons maintenant à la distribution des valeurs propres dans les 
dilatés d'un profil conique de la forme $\Gamma(0,g)\Subset\Omega,$
$\Gamma_{\theta_{1},\theta_{2}}(g,h)$ a été introduit en (\ref{i.26}).
On peut supposer sans perte de généralité que $\inf_{\theta\in[\theta_{1},\theta_{2}]}g(\theta)=1.$

On procède à un  découpage dyadique de $\Gamma(0,\lambda g )$
pour de grandes valeurs de $\lambda.$
Introduisons $k_{0}$ l'entier pour lequel 
$2^{k_{0}}\le\lambda <2^{k_{0}+1}.$ On trouve
\begin{align}\label{w.35}
\Gamma(0,\lambda g)=&
\Gamma(0,1)\cup\left(\bigcup_{k=0}^{k_{0}-1}\Gamma(2^{k},2^{k+1})\right)
\cup \Gamma(2^{k_{0}},\lambda g)\nonumber\\
=&
\Gamma(0,1)\cup\left(\bigcup_{k=0}^{k_{0}-1}2^{k}\Gamma(1,2)\right)
\cup 2^{k_{0}}\Gamma(1,\lambda g/2^{k_{0}}).
\end{align}

\begin{lemme}\label{w12} Supposons $m-\alpha_{1}-\rho-3/4>0.$
Il existe alors $C>0$ tel que quelque soit 
$\tilde{\epsilon}>0$, il existe $k(\tilde{\epsilon})\in \N$ tel que avec une probabilité 
$\ge 1-\tilde{\epsilon}$, on ait
\begin{align}
&\forall \lambda\ge  2^{k(\tilde{\epsilon})},\nonumber\\
&|N(P-Q_{\omega}, \Gamma (2^{k(\tilde{\epsilon})},\lambda g))
-\frac{1}{2\pi}\iint m_{\Gamma (2^{k(\tilde{\epsilon})},\lambda g)}\,dxd\xi|
\le C \lambda^{1/(2m)}\sqrt{\ln\lambda}. \label{w.38}
\end{align}
\end{lemme}
\textbf{Preuve.}
Nous tirons de la section précédente :
avec une probabilité  
$\ge 1-C k^{-\frac{1}{2}}2^{-2\frac{k\gamma_{1}}{m}}$ 
nous avons
\begin{equation}\label{w.39}
|N(P-Q_{\omega},\Gamma(2^{k},2^{k+1}))-\frac{1}{2\pi } \iint m_{\Gamma(2^{k},2^{k+1})}\, dxd\xi | \le C
2^{k/(2m)}\sqrt{k}
\end{equation}
Similairement, 
avec une probabilité  
$\ge 1-C k_{0}^{-\frac{1}{2}}2^{-2\frac{k_{0}\gamma_{2}}{m}}$ 
nous avons pour tout $2^{k_{0}}\le\lambda<2^{k_{0}+1},$
\begin{equation}\label{w.40}
|N(P-Q_{\omega},\Gamma(2^{k_{0}},\lambda g))-\frac{1}{2\pi } \iint m_{\Gamma(2^{k_{0}},\lambda g)}\, dxd\xi | \le C
2^{k_{0}/(2m)}\sqrt{k_{0}}
\end{equation}
(la famille de domaines $\Gamma(1,\lambda g/2^{k_{0}})$ indexée sur $\lambda$  
entre bien dans le cadre du théorème \ref{ss9}.).

Soit $A_{k}$ l'évènement (\ref{w.39}) et $B_{k_{0}}$ l'évènement
(\ref{w.40}).  Pour tout $\gamma_{1},\gamma_{2}$ dans $(0,m-\alpha_{1}-\rho-3/4),$
Nous avons 
\begin{equation}\label{w.41}
\sum_{k=1}^{\infty}
\mathbb{P}(\complement A_{k})+\mathbb{P}(\complement B_{k})=
C\sum_{k=1}^{\infty} k^{-\frac{1}{2}}(2^{-2\frac{k\gamma_{1}}{m}}+
2^{-2\frac{k\gamma_{2}}{m}}) <+\infty.
\end{equation}
Puisque la somme est finie, il existe $k(\tilde{\epsilon})>0$ tel que 
\[\mathbb{P}(\bigcup_{k=k(\tilde{\epsilon})}^{\infty} \complement A_{k} \cup
\complement B_{k})\le \sum_{k=k(\tilde{\epsilon})}^{\infty}
\mathbb{P}(\complement A_{k})+\mathbb{P}(\complement B_{k})<\tilde{\epsilon},\]
impliquant
\[\mathbb{P}(\bigcap_{k=k(\tilde{\epsilon})}^{\infty} (A_{k} \cap
B_{k}))\ge 1-\tilde{\epsilon}.\]

Utilisant le fait que 
\[\sum_{k=k(\tilde{\epsilon})}^{k_{0}}2^{\frac{k}{2m}}\sqrt{\ln 2^{k}}
=\mathcal{O}(1)\lambda^{\frac{1}{2m}}\sqrt{\ln\lambda},
\]
on conclut qu'avec une probabilité $>1-\tilde{\epsilon}$ nous avons
(\ref{w.38}).
\hfill $\square$ \medskip
\begin{theo}\label{w13} Supposons $m_{1}>0.$
Il existe $C_{1}>0$ tel que $\forall\epsilon>0$, il existe 
$M_{\epsilon}\subset\mathcal{M}$ tel que $\mathbb{P}(M_{\epsilon})\ge 1-\epsilon$
et $\forall\omega\in M_{\epsilon},$ il existe $C(\epsilon,\omega)<\infty$ tel que l' on ait
\begin{align*}
&\forall \lambda\ge 0,\\
&| N(P-Q_{\omega}, \lambda\Gamma(0,g)
-\frac{1}{2\pi}\iint m_{\lambda\Gamma(0,g)}|
\le C(\epsilon,\omega)+C_{1} \lambda^{1/(2m)}\sqrt{\ln \lambda}. 
\end{align*}
\end{theo}

Finalement, en prenant $\widetilde{M}=\bigcup_{\epsilon}M_{\epsilon}$, nous avons 
$\mathbb{P}(\tilde{M})=1$
et le 
\begin{corol}\label{w14}
Supposons $m_{1}>0.$
Il existe $C_{1}>0$ et  
$\widetilde{M}\subset\mathcal{M}$ avec  $\mathbb{P}(\widetilde{M})= 1$ tels que
 $\forall\omega\in \widetilde{M}$ on ait
\begin{align*}
&\forall \lambda\ge 0,\\
&|N(PQ_{\omega},\lambda\Gamma(0,g))
-\frac{1}{2\pi}\iint m_{\lambda\Gamma(0,g)}|
\le C(\omega)+C_{1} \lambda^{1/(2m)}\sqrt{\ln \lambda}. 
\end{align*}
\end{corol}

\begin{remarque}\textnormal{
La somme  (\ref{w.39}) est finie, c'est la condition nécessaire pour appliquer le lemme de Borel-Cantelli, rappellé ici
\[\sum_{n}\mathbb{P}(\complement A_{n})<\infty\Rightarrow \mathbb{P}(\lim \inf A_{n})=1,\]
où
\[\omega\in \lim\inf A_{n}\iff \exists\, n(\omega)\in \N,\,\forall k\ge n(\omega),\,
\omega\in A_{k}.\]
Le lemme \ref{w12} et le théorème \ref{w13} n'est autre qu'une redémonstration 
du lemme de Borel-Cantelli appliqué aux évènements $A_{k}$ et $B_{k}$.}
\end{remarque}

\end {document}